\documentclass[11pt]{article}\textwidth 160mm\textheight 235mm
\usepackage{amsfonts}
\usepackage{amssymb}
\usepackage{graphicx}
\usepackage{amsmath}
\usepackage{amsthm}
\usepackage{enumitem}
\usepackage{tikz}
\usepackage{cancel}
\usepackage{todonotes}
\usepackage[normalem]{ulem}
\usetikzlibrary{matrix}
\usetikzlibrary{arrows,shapes}
\usepackage{cite}
\usepackage{hyperref}
\hypersetup{colorlinks=true, urlcolor=blue}
\expandafter\let\expandafter\oldproof\csname\string\proof\endcsname
\let\oldendproof\endproof
\renewenvironment{proof}[1][\proofname]{%
\oldproof[\ttfamily \scshape \bf #1. ]%
}{\oldendproof}

\def\ve{\varepsilon}
\def\epsilon{\varepsilon}
\def\T{{\rm T}}

\def\tilde{\widetilde}
\def\emp{\emptyset}

\def\dom{{\rm dom}\,}
\def\span{{\rm span}\,}
\def\epi{{\rm epi\,}}
\def\rge{{\rm rge\,}}

\def\E{{\cal E}}

\def\N{{\cal N}}

\def\tet{{\tilde{\eta}}}

\def\Lm{{\Lambda}}
\def\O{{\cal O}}

\def\sub{\partial}
\def\B{\mathbb B}

\def\ox{\overline{x}}
\def\oy{\overline{y}}
\def\oz{\overline{z}}
\def\Sm{{ {\cal S}_+^m}}
\def\cl{{\rm cl}\,}

\def\disp{\displaystyle}

\def\H{Hadamard}

\def\tto{\;{\lower 1pt \hbox{$\rightarrow$}}\kern-10pt
\hbox{\raise 2pt\hbox{$\rightarrow$}}\;}
\def\Hat{\widehat}

\def\Bar{\overline}
\def\ra{\rangle}
\def\la{\langle}
\def\ve{\varepsilon}
\def\B{I\!\!B}

\def\R{\mathbb{R}}

\def\X{\mathbb{X}}
\def\Y{\mathbb{Y}}
\def\E{\mathbb{E}}
\def\Rm{\mathbb{R}^{m-1}}
\def\N{I\!\!N}
\def\ox{\bar{x}}
\def\oy{\bar{y}}
\def\oz{\bar{z}}

\def\omu{\bar{\mu}}

\def\ri{\mbox{\rm ri}\,}
\def\tr{\mbox{\rm tr}\,}
\def\int{\mbox{\rm int}\,}
\def\gph{\mbox{\rm gph}\,}
\def\epi{\mbox{\rm epi}\,}
\def\dim{\mbox{\rm dim}\,}
\def\dom{\mbox{\rm dom}\,}
\def\bd{\mbox{\rm bd}\,}
\def\ker{\mbox{\rm ker}\,}

\def\diag{\mbox{\rm diag}\,}

\def\aff{\mbox{\rm aff}\,}

\def\cl{\mbox{\rm cl}\,}
\def\rank{\mbox{\rm rank}\,}

\def\dn{\downarrow}
\def\O{\Omega}
\def\H{{\cal H}}
\def\ph{\varphi}
\def\emp{\emptyset}
\def\st{\stackrel}

\def\lm{\lambda}
\def\olm{\bar\lambda}

\def\dd{\delta}
\def\al{\alpha}
\def\Th{\Theta}
\def\N{I\!\!N}

\def\sce{\setcounter{equation}{0}}
\def\Q{{\cal Q}}

\def\ss{\scriptsize }
\begin{document}
\vspace*{0.5in}
\begin{center}
{\bf CRITICALITY OF LAGRANGE MULTIPLIERS IN VARIATIONAL SYSTEMS}\\[1ex]
BORIS S. MORDUKHOVICH\footnote{Department of Mathematics, Wayne State University, Detroit, MI 48202, USA (boris@math.wayne.edu). Research of this author was partly supported by the National Science Foundation under grants DMS-1512846 and DMS-1808978, and by the Air Force Office of Scientific Research under grant \#15RT0462} and M. EBRAHIM SARABI\footnote{Department of Mathematics, Miami University, Oxford, OH 45065, USA (sarabim@miamioh.edu).}
\end{center}
\vspace*{0.05in}
\small{\bf Abstract.} The paper concerns the study of criticality of Lagrange multipliers in variational systems that has been recognized in both theoretical and numerical aspects of optimization and variational analysis. In contrast to the previous developments dealing with polyhedral KKT systems and the like, we now focus on general nonpolyhedral systems that are associated, in particular, with problems of conic programming. Developing a novel approach, which is mainly based on advanced techniques and tools of second-order variational analysis and generalized differentiation, allows us to overcome principal challenges of nonpolyhedrality and to establish complete characterizations on noncritical multipliers in such settings. The obtained results are illustrated by examples from semidefinite programming.\\[1ex]
{\bf Abbreviated title.} Criticality of multipliers\vspace*{-0.2in}

\newtheorem{Theorem}{Theorem}[section]
\newtheorem{Proposition}[Theorem]{Proposition}
\newtheorem{Remark}[Theorem]{Remark}
\newtheorem{Lemma}[Theorem]{Lemma}
\newtheorem{Corollary}[Theorem]{Corollary}
\newtheorem{Definition}[Theorem]{Definition}
\newtheorem{Example}[Theorem]{Example}
\renewcommand{\theequation}{{\thesection}.\arabic{equation}}
\renewcommand{\thefootnote}{\fnsymbol{footnote}}

\normalsize
\section{Introduction}\sce \label{intro}\vspace*{-0.1in}

This paper is devoted to investigating some core issues of optimization and variational analysis that revolve around {\em criticality} of dual elements (Lagrange multipliers) in the corresponding Karush-Kuhn-Tucker (KKT) systems. The motivation to study multiplier criticality came from applications to convergence rates of primal-dual algorithms of numerical optimization. Then it has been realized that understanding these issues requires a careful theoretical investigation that reveals, in particular, deep interrelations between criticality and other fundamental concepts of variational analysis and generalized differentiation, which are of their own interest.

The notion of criticality, i.e., critical and noncritical Lagrange multipliers, was introduced by Izmailov \cite{iz05} for ${\cal C}^2$-smooth problems of nonlinear programming (NLPs) with equality constraints. It has been recognized from the very beginning that the existence of critical multipliers is the main reason to prevent superlinear convergence of primal iterations in Newtonian methods, since such multipliers persistently attract convergence of dual components. Theoretical and computational issues concerning this phenomenon in nonlinear programs and related variational inequalities were analyzed in many publications and reflected in the monograph by Izmailov and Solodov \cite{is14}. We also refer the reader to their excellent survey \cite{is15}, which is specially devoted to various aspects of multiplier criticality in major primal-dial methods of nonlinear programming; see also the comments by Fischer, Martinez, Mordukhovich, and Robinson to this survey.

A striking property of noncritical Lagrange multipliers is that they yield a certain stability (calmness) property of solution maps to canonical  perturbed KKT systems, which in turn helps to establish {\em superlinear} convergence for Newtonian methods. For instance, Izmailov and Solodov \cite{is12} prove in this way that, in the NLP framework, convergence to a noncritical Lagrange multiplier ensure superlinear rate of convergence of primal-dual iterations in the stabilized sequential quadratic programming (sSQP) method even when the problem is degenerate, i.e., the corresponding set of Lagrange multipliers is not a singleton.

Our recent paper \cite{ms17} conducts a systematic study of criticality for {\em polyhedral} variational systems (generalized KKT) that cover a significantly larger territory than NLPs. Employing advanced tools of second-order variational analysis and generalized differentiation, we obtain therein several characterizations of critical and noncritical multipliers and establish their connections with other fundamental as well as novel  properties of variational systems. In particular, it is shown in \cite{ms17} that the well-recognized and comprehensively characterized property of {\em full stability} of local minimizers in polyhedral problems of constrained optimization allows us to exclude the appearance of critical multipliers associated with such minimizers.

The current paper addresses the study of criticality for the following class of {\em nonpolyhedral} variational systems described in the generalized KKT form
\begin{equation}\label{VS}
\Psi(x,\lm):=f(x)+\nabla\Phi(x)^*\lm=0,\quad\;\lm\in N_\Th\big(\Phi(x)\big),
\end{equation}
where $f\colon\X\to\X$ is a differentiable mapping while $\Phi\colon\X\to\Y$ is a twice differentiable mapping between finite-dimensional spaces, where $\Th\subset Y$ is a closed set with $N_\Th$ standing for its (limiting) normal cone \eqref{2.4}, and where the symbol $^*$ signifies the matrix transposition/adjoint operator. A major source for the generalized KKT system \eqref{VS} comes from the first-order necessary optimality conditions for constrained optimization problems. Indeed, consider a differentiable function $\ph_0:\X\to\R$ and define a constrained optimization problem by
\begin{equation}\label{coop}
\mbox{minimize }\;\ph_0(x)\;\mbox{ subject to }\;\Phi(x)\in\Th,
\end{equation}
where $\Phi$ and $\Th$ are taken from \eqref{VS}. It is well known that system \eqref{VS} with $f:=\nabla\ph_0$ gives us, under a certain constraint qualification, necessary optimality conditions for \eqref{coop}.

Despite a good understanding of noncriticality for systems \eqref{VS} with polyhedral sets $\Th$, not much has been done in the case of nonpolyhedrality. The results established recently in \cite[Theorem~3.3]{zz} and \cite[Proposition~4.2]{lp18} do not provide a satisfactory picture in this regard. Indeed, the assumptions imposed therein are so strong that they may not be satisfied even for classical problems of nonlinear programming.\vspace*{0.02in}

This paper aims at developing a novel approach to the study of critical and noncritical Lagrange multipliers associated with \eqref{VS}, where $\Th$ belongs to a rather general class of regular sets that includes, in particular, all the convex ones. The new notion of {\em semi-isolated calmness} is crucial for our characterizations of noncritical multiplies and subsequent applications. Prior to a detailed consideration of this property, let us emphasize the following: (1) it is strictly weaker than the isolated calmness used, e.g., in \cite{bo94,is14} to justify superlinear convergence of the sequential quadratic programming (SQP) method for nonlinear programs, and (2) it allows us to deal with optimization problems admitting nonunique Lagrange multipliers.

It is important to realize that the generalized KKT systems \eqref{VS} with nonpolyhedral sets $\Th$ fail to satisfy some properties that are granted under polyhedrality. In particular, the semi-isolated calmness property for polyhedral systems \eqref{VS} follows from the uniqueness and noncriticality of Lagrange multipliers. However, it is not the case for nonpolyhedral systems as revealed by Example~\ref{ex3} below. This occurs due to the lack of a certain error bound, which is guaranteed by the Hoffman lemma in polyhedral settings. To overcome this challenge, we first establish new characterizations of {\em uniqueness} of Lagrange multipliers combined with some error bound. This plays a significant role in deriving our main result, Theorem~\ref{uplip}, which provides a complete {\em characterization of noncriticality} under a general reducibility assumption.\vspace*{0.02in}

The rest of the paper is organized as follows. Section~\ref{sect2} recalls some basic concepts of variational analysis and generalized differentiation utilized below. In Section~\ref{sect3} we define critical and noncritical multipliers for system \eqref{VS} together with an extended notion of ${\cal C}^2$-reducibility of $\Th$ and then provide elaborations of these notions for major models of conic programming. Section~\ref{sect4} establishes new characterization of uniqueness of Lagrange multipliers in nonpolyhedral systems. In Section~\ref{sect5} we develop a reduction approach for the study of criticality of multipliers in \eqref{VS} under the ${\cal C}^2$-reducibility of $\Th$ and establish in this way verifiable characterizations of noncritical multipliers with relationships to semi-isolated calmness. Furthermore, we show that the assumptions required for the obtained characterizations are fulfilled under the well-known strict complementarity condition.\vspace*{0.02in}

Our notation and terminology are standard in variational analysis and generalized differentiation; see, e.g., \cite{m06,rw}. Recall that, given a nonempty set $\O$ in $\X$, the notation $\bd\O$, $\int\O$, $\ri\O$, $\cl\O$, $\O^*$, $\aff\O$, and $\span\O$ stands for the boundary, interior,  relative interior, closure, polar, affine hull of $\O$, and the smallest linear subspace containing $\O$, respectively. The symbol $x\st{\O}{\to}\ox$ indicates that $x\to\ox$ with $x\in\O$. By $\B$ we denote the closed unit ball in the space in question while $\B_r(x):=x+r\B$ stands for the closed ball centered at $x$ with radius $r>0$. The indicator function of $\O$ is defined by $\dd_\O(x):=0$ for $x\in\O$ and by $\dd_\O(x):=\infty$ otherwise. Denote by $\diag(a_1,\ldots,a_m)$ an $m\times m$ diagonal matrix whose diagonal entries are $a_1,\ldots,a_m$. We write $x=o(t)$ with $x\in\X$ and $t\in\R_+$ to indicate as usual that $\|x\|/t\to 0$ as $t\dn 0$. Finally, denote by $\R_+$ (respectively, $\R_-$) the set of nonnegative (respectively, nonpositive) real numbers.\vspace*{-0.2in}

\section{Preliminaries from Variational Analysis}
\sce\label{sect2}\vspace*{-0.1in}

In this section we first briefly review, following mainly the books \cite{m06,rw}, basic constructions of variational analysis and generalized differentiation employed in the paper.

Given a set $\O\subset\X$, the (Bouligand-Severi) {\em tangent cone} $T_\O(\ox)$ to $\O$ at $\ox\in\O$ is defined by
\begin{eqnarray}\label{2.5}
T_\Omega(\ox):=\big\{w\in\X\;\big|\;\exists\,t_k{\downarrow}0,\;w_k\to w\;\mbox{ as }\;k\to\infty\;\mbox{with}\;\ox+t_kw_k\in\O\big\}.
\end{eqnarray}
The (Fr\'echet) {\em regular normal cone} $\O$ at $\ox\in\O$ is
\begin{equation}\label{rnc}
\Hat N_\O(\ox):=\big\{v\in\X\;\big|\;\limsup_{x\st{\O}{\to}\ox}\frac{\la v,x-\ox\ra}{\|x-\ox\|}\le 0\big\},
\end{equation}
which can be equivalently described as $\Hat N_\O(\ox)=\T_\O(\ox)^*$. The (limiting/Mordukhovich) {\em normal cone} to $\O$ at $\ox\in\O$ is defined by
\begin{eqnarray}\label{2.4}
N_\O(\ox)=\big\{v\in\X\;\big|\;\exists\,x_k{\to}\ox,\;v_k\to v\;\mbox{ with }\;v_k\in\Hat N_\O(x_k)\big\}.
\end{eqnarray}
If $\O$ is convex, both constructions \eqref{rnc} and \eqref{2.4} reduce to the classical normal cone of convex analysis. The set $\O$ is called (normally) {\em regular} at $\ox\in\O$ if $\Hat N_\O(\ox)=N_\O(\ox)$. In contrast to \eqref{rnc}, the normal cone \eqref{2.4} and the associated constructions for functions and mappings enjoy comprehensive calculus rules based on variational/extremal principles of variational analysis.

Given an extended-real-valued function $f\colon\X\to\Bar\R:=(-\infty,\infty]$ finite at $\ox$, the {\em subdifferential} of $f$ at $\ox$ is defined via the normal cone to its epigraph $\epi f:=\{(x,\al)\in\X\times\R\;|\;\al\ge f(x)\}$ by
\begin{eqnarray}\label{2.6}
\partial f(\ox):=\big\{v\in\X\;\big|\;(v,-1)\in N_{{\scriptsize\epi f}}\big(\ox,f(\ox)\big)\big\}.
\end{eqnarray}
Considering next a set-valued mapping $F\colon\X\tto\Y$ with its domain and graph given by
$$
\dom F:=\big\{x\in\X\;\big|\;F(x)\ne\emp\big\}\;\mbox{ and }\;\gph F:=\big\{(x,y)\in\X\times\Y\;\big|\;x\in F(x)\big\},
$$
define the {\em graphical derivative} of $F$ at $(\ox,\oy)\in\gph F$  by
\begin{equation}\label{gder}
DF(\ox,\oy)(u):=\big\{v\in\Y\;\big|\;(u,v)\in T_{\scriptsize{\gph F}}(\ox,\oy)\big\},\quad u\in\X.
\end{equation}

Finally in this section, we recall the well-posedness properties of set-valued mappings used in what follows. The mapping $F\colon X\tto Y$ is {\em metrically regular} around $(\ox,\oy)\in\gph F$ if there exist $\ell\ge 0$ together with neighborhoods $U$ of $\ox$ and $V$ of $\oy$ such that
\begin{equation}\label{me2}
d(x;F^{-1}(y))\le\ell\,d\big(y;F(x)\big)\;\mbox{ for all }\;(x,y)\in U\times V,
\end{equation}
where $d(x;\O)$ stands for the distance between $x$ and the set $\O$. The {\em metric subregularity} of $F$ at $(\ox,\oy)$ corresponds to the validity of \eqref{me2} with the fixed point $y=\oy$. We say that $F$ is {\em strongly metrically subregular} at $(\ox,\oy)$ if there are $\ell\ge 0$ and a neighborhood $U$ of $\ox$ for which
\begin{equation*}
\|x-\ox\|\le\ell\,d\big(y;F(x)\big)\;\mbox{ whenever }\;x\in U.
\end{equation*}
$F\colon\X\tto\Y$ is {\em calm} at $(\ox,\oy)\in\gph F$ if there are $\ell\ge 0$ and a neighborhood $U$ of $\ox$ such that
\begin{equation}\label{calm-def}
F(x)\cap V\subset F(\ox)+\ell\|x-\ox\|\B\;\mbox{ for all }\;x\in U.
\end{equation}
The {\em isolated calmness} property of $F$ at $(\ox,\oy)$ is defined by
\begin{equation*}
F(x)\cap V\subset\big\{\oy\big\}+\ell\|x-\ox\|\B\;\mbox{ for all}\;x\in U
\end{equation*}
with some $\ell\ge 0$ and a neighborhood $U$ of $\ox$. It is well known that the calmness and isolated calmness of $F$ at $(\ox,\oy)$ are equivalent to the metric subregularity and strong metric subregularity of the inverse mapping $F^{-1}$ at $(\oy,\ox)$, respectively.\vspace*{-0.2in}

\section{Criticality and Reducibility}
\sce\label{sect3}\vspace*{-0.1in}

In this section we first define critical and noncritical multipliers associated with stationary solutions to variational systems of type \eqref{VS}. Then we discuss a modified notion of set reducibility under which criticality can be efficiently investigated in the framework of conic programming.

Given a point $\ox\in\X$ satisfying the {\em stationary condition}
\begin{equation}\label{stat}
0\in f(\ox)+\partial(\dd_\Th\circ\Phi)(\ox),
\end{equation}
we define the set of {\em Lagrange multipliers} associated with $\ox$ by
\begin{equation}\label{laset}
\Lambda(\ox):=\big\{\lm\in\Y\;\big|\;\Psi(\ox,\lm)=0,\;\lm\in N_\Th\big(\Phi(\ox)\big)\big\}.
\end{equation}
Suppose in what follows that $\Lm(\ox)\ne\emp$, which is ensured by a variety of constraint qualification conditions for the system $\Phi(x)\in\Th$ including the metric subregularity of the set-valued constraint mapping $x\mapsto\Phi(x)-\Th$ at $(\ox,0)$.

The following notions of criticality for \eqref{VS} are taken from \cite[Definition~3.1]{ms17}.\vspace*{-0.07in}

\begin{Definition}[\bf critical and noncritical multipliers]\label{crit} Let $\ox$ satisfy the stationery condition \eqref{stat}. The multiplier $\olm\in\Lambda(\ox)$ is {\sc critical} for \eqref{VS} if there is $\xi\in\X$ with $\xi\ne 0$ satisfying
\begin{equation}\label{crc}
0\in\nabla_{x}\Psi(\ox,\olm)\xi+\nabla\Phi(\ox)^*DN_\Th(\Phi(\ox),\olm)\big(\nabla\Phi(\ox)\xi\big).
\end{equation}
The Lagrange multiplier $\olm\in\Lambda(\ox)$ is {\sc noncritical} for \eqref{VS} when the generalized equation \eqref{crc} admits only the trivial solution $\xi=0$.
\end{Definition}\vspace*{-0.05in}
We can reformulate Definition~\ref{crit} via the mapping $G\colon\X\times\Y\tto\X\times \Y$ given by
\begin{equation}\label{GKKT}
G(x,\lm):=\left[\begin{array}{c}
\Psi(x,\lm)\\
-\Phi(x)
\end{array}
\right]+\left[\begin{array}{c}
0\\
N_\Th^{-1}(\lm)
\end{array}
\right].
\end{equation}
It follows from \cite[Theorem~7.1]{ms17} that $\olm\in\Lambda(\ox)$ is noncritical if and only if
\begin{equation}\label{cri2}
(0,0)\in DG\big((\ox,\olm),(0,0)\big)\big(\xi,\eta)\Longrightarrow\xi=0\;\mbox{ for }\;(\xi,\eta)\in\X\times\Y.
\end{equation}
Observe that the stronger implication
\begin{equation*}
(0,0)\in DG\big((\ox,\olm),(0,0)\big)\big(\xi,\eta)\Longrightarrow(\xi,\eta)=(0,0)\;\mbox{ for }\;(\xi,\eta)\in\X\times\Y
\end{equation*}
ensures the property of strong metric subregularity for the mapping $G$ at $\big((\ox,\olm),(0,0)\big)$; see \cite[Theorem~7.1]{ms17} for more details and discussion.\vspace*{0.02in}

The following property of the set $\Th$ in \eqref{VS} is crucial for our subsequent analysis.\vspace*{-0.07in}

\begin{Definition}[\bf reducible sets]\label{defr}
A closed set $\Theta\subset\Y$ is said to be ${\cal C}^2$-{\sc cone reducible} at $\oz=\Phi(\ox)\in\Theta$ to a closed convex subcone $C\subset\E$ of a finite-dimensional space $E$ if there exist a neighborhood ${\cal O}\subset\Y$ of $\oz$ and a ${\cal C}^2$-smooth mapping $h\colon\Y\to\E$ such that
\begin{equation}\label{red}
\Theta\cap{\cal O}=\big\{z\in{\cal O}\;\big|\;h(z)\in C\big\},\quad h(\oz)=0,\;\;\mbox{and}\;\;\nabla h(\oz)\;\mbox{ is surjective}.
\end{equation}
If this holds for all $\oz\in\Theta$, then we say that $\Theta$ is ${\cal C}^2$-cone reducible.
\end{Definition}\vspace*{-0.05in}

Let us discuss this notion and its comparison with the known one in more details.\vspace*{-0.05in}

\begin{Remark}[\bf discussion on reducible sets]{\rm The conventional notion of reducibility from \cite[Definition~3.135]{bs} requires that the convex cone $C$ be pointed. The approach in this paper based on Definition~\ref{defr} does not need this assumption. Moreover, in contrast to \cite[Definition~3.135]{bs} we do not assume that the set $\Th$ is convex; however, \eqref{red} implies that $\Th$ is regular at any $z\in{\cal O}$. Another important point about reducible sets is the requirement that $h(\oz)=0$. This assumption plays a significant role in what follows and cannot be dropped. It helps to reduce our analysis at $\oz$ in $\Th$ to that at $h(\oz)=0$ in another convex cone $C$. Since $N_C(h(\oz))=C^*$, the required inclusion holds automatically for $C$. Thus our approach is to reduce the consideration to $C$, prove the claimed results for this cone, and then return to $\Th$.}
\end{Remark}\vspace*{-0.05in}

The ${\cal C}^2$-cone reducibility of $\Th$ allows us to deduce from the conventional first-order chain rules of variational analysis that for any  $z\in\Th\cap{\cal O}$ with ${\cal O}$ taken from \eqref{red} we have the normal and tangent cone representations
\begin{equation}\label{rch}
N_\Th(z)=\nabla h(z)^*N_C\big(h(z)\big)\quad\mbox{and}\quad T_\Th(z)=\big\{v\in\Y\;\big|\;\nabla h(z)v\in T_C\big(h(z)\big)\big\}.
\end{equation}
Let us now consider in more details the three important cases of the variational system \eqref{VS} where $\Th$ therein is one of the following sets:\vspace*{-0.1in}
\begin{itemize}[noitemsep]
\item convex polyhedral set;
\item the second-order cone;
\item the cone of positive semidefinite symmetric matrices.
\end{itemize}\vspace*{-0.1in}
It is well known that these sets are ${\cal C}^2$-cone reducible; see \cite[Examples~3.139 and 3.140]{bs}. Below we provide simplified and constructive proofs for these reductions. Our first example concern polyhedral sets, where--in contrast to \cite[Examples~3.139]{bs}--we explicitly construct
$h$ in \eqref{red} as an affine mapping, which is used in our subsequent analysis.\vspace*{-0.07in}

\begin{Example}[\bf convex polyhedra]\label{poly}{\rm Let $\Y=\R^m$, and let $\Th$ in \eqref{VS} be a convex polyhedral set with $\oz\in\Th$. We intend to show that $\Th$ is ${\cal C}^2$-cone reducible at $\oz$. Denote $s:=\dim\span\{N_\Th(\oz)\}$ and let $A$ be the matrix of linear isometry from $\mathbb{R}^m$ into $\mathbb{R}^s\times\mathbb{R}^{m-s}$ under which $A^*(\span\{N_\Th(\oz)\})=\mathbb{R}^s\times\{0\}$. Represent $\oy:=A^{-1}\oz$ as $\oy=(\oy_s,\oy_{m-s})\in\R^s\times\R^{m-s}$ and define the set $D\subset\R^s$ by
\begin{equation*}
D:=\big\{x\in\R^s\;\big|\;A(x,\oy_{m-s})\in\Th\big\},
\end{equation*}
which is clearly a convex polyhedron. Construct now an $s\times m$ matrix $B$ by deleting the last $m-s$ rows of the $m\times m$ matrix $A^{-1}$.
Using the same arguments as \cite[Lemma~3.2]{ms16} gives us
\begin{equation*}
\Th=\big\{z\in\R^m\;\big|\;Bz\in D\big\}.
\end{equation*}
Since $D$ is convex polyhedron, it follows from \cite[Theorem~2E.3]{dr} that there is a neighborhood $U$ of $0\in\R^s$ for which $T_D(B\oz)\cap U=\big(D-B\oz\big)\cap U$. Define further $h(z):=Bz-B\oz$ for any $z\in\R^m$ and find by the continuity of $h$ a neighborhood $O$ of $\oz$ such that $h(z)=h(z)-h(\oz)\in U$ whenever $z\in O$. Combining all the above tells us that
\begin{equation*}
\Th\cap O=\big\{z\in O\;\big|\;h(z)\in C\big\}\;\mbox{ with }\;C:=T_D(B\oz).
\end{equation*}
It is easy to check that the constructed mapping $h$ and the convex cone $C\subset\R^s$ with $s=\dim\span\{N_\Th(\oz)\}$ satisfy \eqref{red}, and thus the set $\Th$ is ${\cal C}^2$-cone reducible.}
\end{Example}\vspace*{-0.05in}

The second example addresses a nonpolyhedral cone, which generates an important class of problems of second-order cone programming (SOCP).\vspace*{-0.07in}

\begin{Example}[\bf second-order cone]\label{socp} {\rm Let $\Y=\R^m$, and let $\Th:=\Q\subset\Y$, where $\Q$ is the {\em second-order/Lorentz/ice-cream cone} defined by
\begin{equation}\label{Qm}
\Q:=\big\{s=(s_r,s_m)\in\R^{m-1}\times\R\;\big|\;\Vert s_r\Vert\le s_m\big\}.
\end{equation}
It follows from \cite[Lemma~15]{br} that the second-order cone $\Q$ is ${\cal C}^2$-cone reducible at $\oz\in\Q$ to
\begin{eqnarray*}
C:=\left\{\begin{array}{ll}\Q&\mbox{if}\;\oz=0,\\
\{0\}&\mbox{if}\;\oz\in(\int\Q)\setminus\{0\},\\
\R_{-}&\mbox{if}\;\oz\in(\bd\Q)\setminus\{0\}.
\end{array}
\right.
\end{eqnarray*}
We represent in what follows an element $y\in\Q$ as $y=(y_r,y_m)$ with $y_m\in\R$ and $y_r\in\Rm$. The reduction mapping
$h$ can be defined as
\begin{equation}\label{eqq2}
h(z):=\begin{cases}
z&\textrm{if }\;\oz=0,\\
0\in\R&\textrm{if }\;\oz\in\int\Q,\\
\|z_r\|^2-z_m^2&\textrm{if }\;\oz\in(\bd\Q)\setminus\{0\}
\end{cases}
\end{equation}
for all vectors $z$ in a neighborhood of $\oz$. Picking $z=(z_r,z_m)\in\Q$ and $\lm=(\lm_r,\lm_m)\in N_\Q(z)$, we construct the matrix $\H(z,\lm)$ by
\begin{equation}\label{eqq1}
\H(z,\lm):=\begin{cases}
-\dfrac{\lambda_m}{z_m}\,\diag(\underbrace{1,\ldots,1}_{m-1\;{\rm times}},-1)&\textrm{if }\;z=(z_r,z_m)\in(\bd\Q)\setminus\{0\},\\
0&\textrm{if }\;z\in\big[(\int\Q)\cup\{0\}\big].
\end{cases}
\end{equation}
This matrix appears as the curvature term of the second-order cone $\Q$ in Proposition~\ref{gdno}.}
\end{Example}\vspace*{-0.05in}

Next we consider a more involved cone $\Th$ is \eqref{VS}, which generates problems of {\em semidefinite programming} (SDP) that are highly important in applications.\vspace*{-0.07in}

\begin{Example}[\bf positive semidefinite cone]\label{sdp}{\rm Let $\Y:={\cal S}^m$ be the space of $m\times m$ symmetric matrices, which is conveniently treated via the inner product
$$
\la A,B\ra:=\tr AB
$$
with $\tr AB$ standing for the sum of the diagonal entries of $AB$. This inner product induces a norm on ${\cal S}^m$ known as the {\em Frobenius/Hilbert-Schmidt norm} and defined by
\begin{equation*}
\|A\|:=\big(\sum_{i,j=1}^{m}a_{ij}^2\big)^\frac{1}{2}\;\mbox{ with }\;A:=(a_{ij}).
\end{equation*}
Given $A,B\in\Sm$, it is not hard to see that $\la A,B\ra=0$ if and only if $AB=0$. For a matrix $A\in{\cal S}^m$, denote by $A^\dagger$ the {\em Moore-Penrose pseudo-inverse} of $A$. In this case we have $\Th=\Sm$, where $\Sm$ is the cone of $m\times m$ positive semidefinite symmetric matrices. Denote $\rank A=:p$ for $A\in\Sm$ and consider the following two cases. In the case where $p=m$ the matrix $A$ is positive definite and hence belongs to the interior of $\Sm$. Then it is easy to observe that $\Sm$ is ${\cal C}^2$-cone reducible at $A$ to $\{0\}$ with the reduction mapping $h\colon\Sm\to\{0\}$ defined by $h(B):=0$ for $B$ in a neighborhood of $A$. In the case where $p<m$ we know from \cite[Example~3.140]{bs} that $\Sm$ is ${\cal C}^2$-cone reducible at $A$ to ${\cal S}_+^{m-p}$ via the mapping $h\colon{\cal S}^m\to{\cal S}^{m-p}$ defined by $h(B):=U(B)^*BU(B)$; see  \cite[Example~3.140]{bs} for the definition of $U(B)$ and more details on this mapping. It follows from \cite[Example~3.98]{bs} that $h(A)=U(A)^*AU(A)=\al I_{m-p}$, where $\al$ is the smallest eigenvalue of $A$ and where $I_p$ stands for the $(m-p)\times(m-p)$ identity matrix. Since $p<m$, we have that $\al=0$ and thus $h(A)=0$, which indeed shows that $h$ satisfies in \eqref{red}.}
\end{Example}\vspace*{-0.05in}

The next result calculates the graphical derivative of the normal cone mapping (which is a primal-dual construction of second-order variational analysis) generated by reducible sets $\Th$. This is instrumental for the study of multiplier criticality in such settings. Recall that the {\em critical cone} to $\Th$ at $z\in\Th$ for $\lm\in N_\Th(z)$ is defined by
\begin{equation}\label{cri3}
K_\Th(z,\lm)=T_\Th(z)\cap\{\lm\}^\perp.
\end{equation}\vspace*{-0.3in}

\begin{Proposition}[\bf graphical derivative of normal cones to reducible sets]\label{gdno}
Let $(\oz,\olm)\in\gph N_\Th$, and let $\Th$ be ${\cal C}^2$-{cone reducible} at $\oz$ to a closed convex cone $C$.
Then the graphical derivative of the normal cone mapping $N_\Th$ is calculated by
\begin{equation}\label{gdf}
DN_\Th(\oz,\olm)(u)=\nabla^2\la\bar\mu,h\ra(\oz)u+N_{\ss{K_\Th(\oz,\olm)}}(u)\;\mbox{ for all }\;u\in\Y
\end{equation}
via the critical cone \eqref{cri3}, where $\bar\mu$ is the unique solution to the system
\begin{equation}\label{redl}
\olm=\nabla h(\oz)^*\bar\mu\;\mbox{ and }\;\bar\mu\in N_C\big(h(\oz)\big),
\end{equation}
and where $h$ is taken from \eqref{red}. If $\Th$ is a convex polyhedron in $\Y=\R^m$, then we have $\nabla^2\la\bar\mu,h\ra(\oz)u=0$ as $u\in\R^m$ for the curvature term in \eqref{gdf}. If $\Th=\Q\subset\Y=\R^m$, then
\begin{equation}\label{curq}
\nabla^2\la\bar\mu,h\ra(\oz)u=\H(\oz,\olm)u\quad\mbox{for all}\;\;u\in\R^m.
\end{equation}
Finally, in the SDP case where $\Y={\cal S}^m$ and $\Th=\Sm$ we have the representation
\begin{equation}\label{curs}
\nabla^2\la\bar\mu,h\ra(\oz)u=-2\olm u\oz^\dagger\quad\mbox{for all}\;\;u\in{\cal S}^m.
\end{equation}
\end{Proposition}\vspace*{-0.1in}
\begin{proof} Since $\olm\in\Lambda(\ox)$ and $\nabla h(\oz)$ is surjective, the normal cone representation in \eqref{rch} implies that there is  a unique vector $\bar\mu\in N_C(h(\oz))$ such that $\olm=\nabla h(\oz)^*\bar\mu$. This allows us to deduce \eqref{gdf} from \cite[Corollary~4.5]{gm17}. To calculate the curvature term for the second-order cone $\Q$, we get from \eqref{eqq2} that $\nabla^2\la\bar\mu,h\ra(\oz)u=0$ if $\oz\in[\int\Q]\cup\{0\}$, which verifies \eqref{curq} in this case due to \eqref{eqq1}. If $\oz\in(\bd\Q)\setminus\{0\}$, it follows from \eqref{eqq2} that
\begin{equation*}
h(y)=\|y_r\|^2-y_m^2\;\mbox{ whenever }\;y=(y_r,y_m)\in\Rm\times\R.
\end{equation*}
Since $\omu\in N_{C}(h(\oz))$ with $C=\R_-$, we get $\omu\in\R_+$ and thus conclude from \eqref{redl} that
\begin{equation*}
\olm=\nabla h(\oz)^*\omu=\omu\begin{pmatrix}
2\oz_r\\
-2\oz_m\\
\end{pmatrix},
\end{equation*}
which in turn yields $\omu=-\dfrac{\olm_m}{2\oz_m}$. On the other hand, the direct calculations lead us to
$$
\nabla^2\la\bar\mu, h\ra(\oz)=\omu\,\diag(\underbrace{2,\ldots,2}_{m-1\;{\rm times}},-2)=-\dfrac{\olm_m}{\oz_m}\,\diag(\underbrace{1,\ldots,1}_{m-1\;{\rm times}},-1).
$$
Using now \eqref{eqq1} gives us \eqref{curq} in the case where $\oz\in(\bd\Q)\setminus\{0\}$. To calculate the curvature term for $\Sm$, we employ \cite[equation~(66)]{br2} and get
$$
\la\nabla^2\la\bar\mu,h\ra(\oz)u,u\ra=-2\la\olm,u\oz^\dagger u\ra\;\mbox{ for all }\;u\in{\cal S}^m.
$$
Differentiating both sides above with respect to $u$ brings us to
$$
\nabla^2\la\bar\mu,h\ra(\oz)u=-\dfrac{\sub\la\olm,u\oz^\dagger u\ra}{\sub u}=-\olm u\oz^\dagger-\oz^\dagger u\olm=-2\olm u\oz^\dagger,
$$
which justifies \eqref{curs} and thus completes the proof of the proposition.
\end{proof}\vspace*{-0.05in}

As an immediate consequence of Definition~\ref{crit} and Proposition~\ref{gdno}, we arrive at the following equivalent description of critical multipliers for \eqref{VS} when $\Th$ is a ${\cal C}^2$-cone reducible set.\vspace*{-0.058557in}

\begin{Corollary}[\bf equivalent description of critical multipliers]\label{crit-eq} Let $\ox$ satisfy the stationery condition \eqref{stat}, let $\olm\in\Lambda(\ox)$, let $\Th$ be ${\cal C}^2$-cone reducible at $\oz:=\Phi(\ox)$ to a closed convex cone $C$, and let $\bar\mu$ be a unique solution to \eqref{redl}. Then $\olm$ is critical for \eqref{VS} if and only if the system
\begin{equation*}
\nabla_{x}\Psi(\ox,\olm)\xi+\nabla\Phi(\ox)^*\eta=0\;\mbox{ and }\;\eta-\nabla^2\la\bar\mu,h\ra(\oz)\nabla\Phi(\ox)\xi\in N_{\ss{K_\Th(\oz,\olm)}}\big(\nabla\Phi(\ox)\xi\big)
\end{equation*}
admits a solution $(\xi,\eta)\in\X\times\Y$ such that $\xi\ne 0$.
\end{Corollary}\vspace*{-0.05in}

As mentioned in Section~\ref{intro}, KKT systems corresponding to problems of constrained optimization \eqref{coop} clearly belong to class \eqref{VS}. The Lagrangian for \eqref{coop} is defined by
\begin{equation*}
L(x,\lm):=\ph_0(x)+\la\Phi(x),\lm\ra,
\end{equation*}
while the set of Lagrange multipliers for \eqref{coop} associated with a feasible solution $\ox$ is given by
\begin{equation*}
\Lambda_{\ss{c}}(\ox):=\big\{\lm\in\Y\;\big|\;\nabla_x L(\ox,\lm)=0,\;\lm\in N_\Th\big(\Phi(\ox)\big)\big\}.
\end{equation*}
Let $(\oz,\olm)\in\gph N_\Th$ with $\oz=\Phi(\ox)$, and let $\Th$ be ${\cal C}^2$-cone reducible at $\oz$ to the closed convex cone $C$.
Given $\olm\in\Lambda_{\ss{c}}(\ox)$, we formulate the {\em second-order sufficient condition} for \eqref{coop} as
\begin{equation}\label{ssoc}
\begin{cases}
\big\la\nabla^2_{xx}L(\ox,\olm)u,u\big\ra+\big\la\nabla^2\la\bar\mu,h\ra(\oz)\nabla\Phi(\ox)u,\nabla\Phi(\ox)u \big\ra>0\\
\;\mbox{for all }\;0\ne u\in\X\;\mbox{ with }\;\nabla\Phi(\ox)u\in K_\Th(\oz,\olm),
\end{cases}
\end{equation}
where $h$ and $\omu$ are taken from \eqref{red} and \eqref{redl}, respectively. When $\Y=\R^m$ and $\Th=\Q$, the curvature term in \eqref{ssoc} is calculated in Proposition~\ref{gdno} as $\la\nabla^2\la\bar\mu,h\ra(\oz)u,u\ra=\la\H(\oz,\olm)u,u\ra$ for all $u\in\Y$. If $\Y={\cal S}^m$ and $\Th=\Sm$, the curvature term in \eqref{ssoc} reduces by Proposition~\ref{gdno} to $\la\nabla^2\la\bar\mu,h\ra(\oz)u,u\ra=-2\la\olm,u\Phi(\ox)^\dagger u\ra$ for all $u\in\Sm$. Note that \eqref{ssoc} can be stronger than the classical second-order sufficient condition for \eqref{coop} given by
\begin{equation*}
\begin{cases}
\disp{\sup_{\olm\in\Lambda_{\ss{c}}(\ox)}\big\{\big\la\nabla^2_{xx}L(\ox,\olm)u,u\big\ra+\big\la\nabla^2 \la\bar\mu,h\ra(\oz)\nabla\Phi(\ox)u,\nabla \Phi(\ox)u\big\ra\big\}>0}\\
\;\mbox{for all }\;0\ne u\in\X\;\mbox{ with }\;\nabla\Phi(\ox)u\in K_\Th(\oz,\olm),
\end{cases}
\end{equation*}
if the set of Lagrange multipliers is not a singleton. However, an advantage of \eqref{ssoc} is that it provides a sufficient condition for noncriticality of Lagrange multipliers. Example~\ref{ex1} confirms that it may be much easier to justify noncriticality by using the second-order sufficient condition \eqref{ssoc} than working with definition \eqref{crc} or its simplification from Corollary~\ref{crit-eq}.\vspace*{-0.07in}

\begin{Proposition}[\bf sufficient condition for noncriticality of a Lagrange multipliers]\label{snoc}
Let $\ox$ be a feasible solution to \eqref{coop}, let $\olm\in\Lambda_{\ss{c}}(\ox)$, and let $\Th$ be ${\cal C}^2$-cone reducible at $\oz=\Phi(\ox)$ to a closed convex cone $C$. If the second-order sufficient condition \eqref{ssoc} holds, then $\ox$ is a strict local minimizer for \eqref{coop} and the Lagrange multiplier $\olm$ is noncritical.
\end{Proposition}\vspace*{-0.15in}
\begin{proof}
The first fact is a well-known result, which follows, e,g., from \cite[Theorem~3.86]{bs}. The noncriticality of $\olm$ under \eqref{ssoc} can be verified  directly while arguing by contradiction.
\end{proof}\vspace*{-0.05in}

Let us now present an SDP example borrowed from Shapiro \cite[Example~4.5]{sh} who constructed it for different purposes. In our case it shows via Proposition~\ref{snoc} that the unique Lagrange multiplier is noncritical.\vspace*{-0.07in}

\begin{Example}[\bf SDP]\label{ex1}{\rm Consider the semidefinite program with $\X=\R^2$, $\Y={\cal S}^2$, and $\Th={\cal S}^2_+$:
\begin{equation}\label{sdp1}
\mbox{minimize}\;x_1+\frac{1}{2}x_1^2+\frac{1}{2}x_2^2\;\mbox{ subject to }\;\Phi(x_1,x_2)\in\Th,
\end{equation}'
where $\Phi\colon\R^2\to\Y$ is defined by $\Phi(x_1,x_2):=\diag(x_1,x_2)$. The feasible set of this problem can be written as
$\{(x_1,x_2)\in\R^2\;|\;x_1\ge 0,\;x_2\ge 0\}$. This shows that $\ox:=(0,0)$ is a unique optimal solution  to \eqref{sdp1}.
Picking $\olm\in\Lambda_{\ss{c}}(\ox)$, we see that $\olm$ satisfies the first-order optimality conditions
$$
\nabla_xL(\ox,\olm)=0,\;\la\olm,\Phi(\ox)\ra=0,\;\mbox{ and }\;\olm\in S_-^2.
$$
They imply that $\olm=\diag(-1,0)$, and so the set of Lagrange multipliers is a singleton. It follows from $\Phi(\ox)=\diag(0,0)$ that
$\Phi(\ox)^\dagger=\diag(0,0)$. Thus
\begin{eqnarray*}
&&\big\la\nabla^2_{xx}L(\ox,\olm) u,u\big\ra+\big\la\nabla^2\la\bar\mu,h\ra(\oz)\nabla\Phi(\ox)u,\nabla\Phi(\ox)u\big\ra =\big\la\nabla^2_{xx}L(\ox,\olm)u,u\big\ra\\
&&+2\big\la\olm,\nabla\Phi(\ox)u\Phi(\ox)^\dagger\nabla\Phi(\ox)u \big\ra=\big\la\diag(1,1)u,u\big\ra=\|u\|^2>0\;\mbox{ for all }\;0\ne u\in\R^2,
\end{eqnarray*}
which verifies that the second-order sufficient condition \eqref{ssoc} holds for $\olm$. Employing now Proposition~\ref{snoc} tells us that the unique Lagrange multiplier $\olm$ is noncritical.}
\end{Example}\vspace*{-0.02in}

When the set $\Th$ is ${\cal C}^2$-cone reducible at $\oz=\Phi(\ox)$ to a closed convex cone $C$, it is useful to consider a counterpart of \eqref{VS} for the closed convex cone $C$ from \eqref{red} written as
\begin{equation}\label{rVS}
\Psi^r(x,\mu):=f(x)+\nabla\big(h\circ\Phi\big)(x)^*\mu=0\;\mbox{ and }\;\mu\in N_C\big((h\circ\Phi)(x)\big)
\end{equation}
with $(x,\mu)\in\X\times\E$. The set of Lagrange multipliers for the {\em reduced variational system} \eqref{rVS} associated with a stationary point $\ox$ from \eqref{stat} is defined by
\begin{equation*}
\Lambda^r(\ox):=\big\{\mu\in\E\;\big|\;\Psi^r(\ox,\mu)=0,\;\mu\in N_C\big((h\circ\Phi)(\ox)\big)\big\}.
\end{equation*}
Since $\nabla h(\oz)$ is surjective, we get the relationship
\begin{equation}\label{lare}
\Lambda(\ox)=\nabla h(\oz)^*\Lambda^r(\ox),
\end{equation}
which is largely exploited below.\vspace*{-0.2in}

\section{Uniqueness and Stability of Lagrange Multipliers}
\sce\label{sect4}\vspace*{-0.1in}

This section is devoted to establishing necessary and sufficient conditions for the uniqueness of Lagrange multipliers in nonpolyhedral systems \eqref{VS} combined with their certain error bound. Besides being of its own interest, this issue is very instrumental for characterizing
noncritical multipliers in the next section. Given a stationary point $\ox$ from \eqref{stat}, define the {\em Lagrange multiplier mapping} $M_{\ox}\colon\X\times\Y\tto\Y$ associated with $\ox$ by
\begin{equation}\label{lagmap}
M_{\ox}(v,w):=\big\{\lm\in\Y\;\big|\;(v,w)\in G(\ox,\lm)\big\}\;\mbox{ for all }\;(v,w)\in\X\times\Y,
\end{equation}
where $G$ is taken from \eqref{GKKT}. It is easy to see that $M_{\ox}(0,0)=\Lambda(\ox)$, where $\Lambda(\ox)$ is the set of Lagrange multipliers at $\ox$ defined in \eqref{laset}.

The following theorem provides characterizations of the uniqueness of Lagrange multipliers in \eqref{VS} together with some error bound and calmness properties, which are automatic for polyhedral systems. In particular, in the case of NLPs the obtained characterizations of uniqueness reduce to the strong Mangasarian-Fromovitz constraint qualification (SMFCQ); see \cite[page~11]{is14} for more details. When $\Y=\R^m$ and the set $\Th$ is the second-order cone $\Q$, a similar result has been recently established in \cite[Theorem~4.5]{hms}. Further discussions are given in Remark~\ref{rem1}. \vspace*{-0.07in}

\begin{Theorem}[\bf characterizations of uniqueness and stability of Lagrange multipliers]\label{unique} Let $\ox$ fulfill the stationery condition \eqref{stat}, let $\Th$ be regular at $\ox$, and let $\olm\in\Lambda(\ox)$. Then we have the following equivalent assertions:\vspace*{-0.1in}
\begin{itemize}[noitemsep]
\item[\bf{(i)}] The Lagrange multiplier $\olm$ is unique and there exist constants $\ell\ge 0$ and $\ve>0$ ensuring the error bound
estimate\vspace{-0.2cm}
\begin{equation}\label{gf01}
d\big(\lm;\Lambda(\ox)\big)\le\ell\;\big(\|\Psi(\ox,\lm)\|+d\big(\Phi(\ox);N_\Th^{-1}(\lm)\big)\big)\;\mbox{ for all }\;\lm\in\B_\ve(\olm).
\end{equation}
\item[\bf{(ii)}] The Lagrange multiplier $\olm$ is unique and the mapping ${M}_{\ox}$ from \eqref{lagmap} is calm at $((0,0),\olm)$.
\item[\bf{(iii)}] The Lagrange multiplier mapping $M_{\ox}$ is isolatedly calm at $((0,0),\olm)$.
\item[\bf{(iv)}] The dual qualification condition is satisfied:
\begin{equation}\label{gf02}
DN_\Th(\Phi(\ox),\olm)(0)\cap\ker\nabla\Phi(\ox)^*=\{0\}.
\end{equation}
\end{itemize}
\end{Theorem}\vspace*{-0.15in}
\begin{proof} Assertions (i) and (ii) are equivalent by the definitions. To proceed further, denote $G_{\ox}(\lm):=G(\ox,\lm)$ and see that $G_{\ox}^{-1}=M_{\ox}$. Then (i) amounts to saying that the mapping $G_{\ox}$ is strongly metrically subregular at $(\olm,(0,0))$.
Indeed, the validity of (i) clearly yields the blue strong subregularity property of $G_{\ox}$ at $(\olm,(0,0))$. Conversely, the latter property tells us that \eqref{gf01} holds and that for some $\ve>0$ we get the equalities
$$
M_{\ox}(0,0)\cap\B_\ve(\olm)=G_{\ox}^{-1}(0,0)\cap\B_\ve(\olm)=\{\olm\}.
$$
It follows from the regularity of $\Th$ at $\ox$ that $M_{\ox}$ is convex-valued. Thus $\Lambda(\ox)=M_{\ox}(0,0)=\{\olm\}$, which gives us (i).
Since $G_{\ox}^{-1}=M_{\ox}$, the strong metric subregularity of $G_{\ox}$ at $(\olm,(0,0))$ means the isolated calmness of $M_{\ox}$ at $((0,0),\olm)$, and therefore we have (i)$\Longleftrightarrow$(iii).

It remains to verify the equivalence between (iii) and (iv). Calculating the graphical derivative of $G_{\ox}$ due to structure \eqref{GKKT} gives us
$$
DG_{\ox}\big(\olm,(0,0)\big)(\eta)=\left[\begin{array}{c}
\nabla\Phi(\ox)^*\eta\\
0
\end{array}
\right]+\left[\begin{array}{c}
0\\
DN_\Th^{-1}\big(\olm,\Phi(\ox)\big)(\eta)
\end{array}
\right]\quad\mbox{for all}\;\;\eta\in\Y.
$$
Since the graph of $G_{\ox}$ is closed, we deduce from \cite[Theorem~4E.1]{dr} that $G_{\ox}$ is strongly metrically subregular at $(\olm,(0,0))$  if and only if the implication
$$
(0,0)\in DG_{\ox}\big(\olm,(0,0)\big)(\eta)\Longrightarrow\eta=0
$$
holds. This amounts to saying that
$$
\eta\in DN_\Th\big(\Phi(\ox),\olm\big)(0)\cap\ker\nabla\Phi(\ox)^*\Longrightarrow\eta=0.
$$
The latter verifies the equivalence between (iii) and (iv), and thus completes the proof.
\end{proof}\vspace*{-0.15in}

\begin{Remark}[\bf discussion on error bounds]\label{rem1}{\rm It can be checked by the direct calculation that in the case of NLPs in \eqref{VS} the dual qualification condition \eqref{gf02} reduces to SMFCQ. In the latter framework the error bound estimate \eqref{gf01} always holds and can be derived by applying the classical Hoffman lemma (see, e.g., \cite[Lemma~3C.4]{dr}) to the Lagrange multiplier mapping $M_{\ox}$ from \eqref{lagmap}.  This explains why for nonlinear programming problems the uniqueness of Lagrange multipliers and  SMFCQ are equivalent. More broadly, if $\Th$ is a convex polyhedral set, we can show that \eqref{gf01} holds automatically. Indeed, we know from convex analysis that $N_\Th^{-1}=\sub\dd^*_\Th$. Thus it follows from \cite[Theorem~11.14]{rw} that $\dd^*_\Th$ is convex piecewise linear in the sense of \cite[Definition~2.47]{rw}), and so its subdifferential mapping is outer/upper Lipschitzian due to Robinson's seminal result \cite{rob1}. This allows us to justify the error bound estimate \eqref{gf01} when $\Th$ is a (convex) {\em polyhedron}. It is not hard to go further and show that if the normal cone $N_\Th$ is replaced by the subdifferential mapping of a {\em convex piecewise linear-quadratic function} from \cite[Definition~10.20]{rw}, then estimate \eqref{gf01} also automatically fulfills.}
\end{Remark}\vspace*{-0.05in}

The result of \cite[Proposition~4.50]{bs} tells us that the {\em strong Robinson constraint qualification} (SRCQ) defined in primal terms by
\begin{equation}\label{sqc}
\nabla\Phi(\ox)\X+T_\Th\big(\Phi(\ox)\big)\cap\{\olm\}^\bot=\Y
\end{equation}
(this terminology was suggested in \cite{dsz}) provides a sufficient condition for the uniqueness of Lagrange multipliers Lagrange in constrained optimization with $\Th$ being a closed, convex while not necessarily ${\cal C}^2$-cone reducible set. On the other hand, the novel dual qualification condition \eqref{gf02} addresses the generalized KKT systems \eqref{VS} that appear in a broader framework than constrained optimization and occurs to be sufficient for the uniqueness of multipliers therein for reducible sets $\Th$. As we have recently proved in \cite[Theorem~4.5]{hms}, both constraint qualifications are equivalent when $\Y=\R^m$ and $\Th$ is the second-order cone $\Q$. Now we extend this result to the general case where $\Th$ is any ${\cal C}^2$-cone reducible set, which may not even be convex. \vspace*{-0.07in}

\begin{Proposition}[\bf equivalence between and dual constraint qualifications under reducibility]\label{pdcq}
Let $\ox$ satisfy the stationery condition \eqref{stat}, let $\olm\in\Lambda(\ox)$, and let $\Th$ be ${\cal C}^2$-cone reducible at $\Phi(\ox)$ to a closed convex cone $C$. Then the dual qualification condition \eqref{gf02} is equivalent to SRCQ \eqref{sqc}.
\end{Proposition}\vspace*{-0.17in}
\begin{proof} It follows from \eqref{gdf} that
\begin{equation}\label{grader}
DN_\Th(\oz,\olm)(0)=N_{\ss{K_\Th(\oz,\olm)}}(0)\;\mbox{ with}\;\oz=\Phi(\ox).
\end{equation}
Assuming the validity of SRCQ, we get the equalities
$$
K_\Th(\oz,\olm)^*\cap\ker\nabla\Phi(\ox)^*=\big(T_\Th(\oz)\cap\{\olm\}^\bot\big)^*\cap\ker\nabla\Phi(\ox)^*=\big( T_\Th(\oz)\cap\{\olm\}^\bot +\nabla\Phi(\ox)\X\big)^*=\{0\}.
$$
Combining this with \eqref{grader} clearly yields \eqref{gf02}. Conversely, assuming \eqref{gf02} and appealing again to \eqref{grader} tell us that
\begin{eqnarray*}
\cl\big(\nabla\Phi(\ox)\X+T_\Th(\oz)\cap\{\olm\}^\bot\big)&=&\big(K_\Th(\oz,\olm)^*\cap\ker\nabla\Phi(\ox)^*\big)^*\\
&=&\big(DN_\Th)(\oz,\olm)(0)\cap\ker\nabla\Phi(\ox)^*\big)^*=\Y.
\end{eqnarray*}
Since the set $\nabla\Phi(\ox)\X+T_\Th(\Phi(\ox))\cap\{\olm\}^\bot$ is convex, it has nonempty relative interior. Hence it follows from \cite[Proposition~2.40]{rw} that the relationships
\begin{eqnarray*}
\Y=\ri(\Y)&=&\ri\big[{\rm cl}\big(\nabla\Phi(\ox)\X+T_\Th\big(\Phi(\ox)\big)\cap\{\olm\}^\bot\big)\big]\\
&=&\ri\big(\nabla\Phi(\ox)\X+T_\Th\big(\Phi(\ox)\big)\cap\{\olm\}^\bot\big)\\
&\subset&\big(\nabla\Phi(\ox)\X+T_\Th\big(\Phi(\ox)\big)\cap\{\olm\}^\bot\big)
\end{eqnarray*}
are satisfied, which therefore completes the proof.
\end{proof}\vspace*{-0.07in}

We highlight here that Theorem~\ref{unique} seems to be the first result in the literature, which provides not only sufficient but also necessary conditions for the uniqueness of Lagrange multipliers in the general framework of \eqref{VS}. As mentioned above, the uniqueness of Lagrange multipliers for NLPs is fully characterized by SMFCQ. However, it follows from Theorem~\ref{unique} that in the general setting of \eqref{VS} the validity of such a result demands that the Lagrange multiplier mapping $M_{\ox}$ be calm. Is the calmness of the latter mapping essential for the validity of Theorem~\ref{unique}? The next example confirms that it is the case, in particular, forb the SDPs.\vspace*{-0.05in}

\begin{Example}[\bf failure of the dual qualification condition for SDPs with unique Lagrange multipliers]\label{ex2}{\rm Consider SDP \eqref{sdp1} from Example~\ref{ex1}, where $\Th={\cal S}_+^2$ is ${\cal C}^2$-cone reducible. To verify that the dual qualification condition \eqref{gf02} fails, observe from \eqref{grader} that
\begin{equation*}
DN_\Th\big(\Phi(\ox),\olm\big)(0)\cap\ker\nabla\Phi(\ox)^*=K_{{\cal S}_+^2}(\oz,\olm)^*\cap\ker\nabla\Phi(\ox)^*,
\end{equation*}
where $\oz:=\Phi(\ox)=\diag(0,0)$ and $\olm=\diag(-1,0)$. We calculate the critical cone $K_{S_+^2}(\oz,\olm)$ by
\begin{equation*}
K_{{\cal S}_+^2}(\oz,\olm)=\big\{u\in{\cal S}_+^;\big|\;\la u,\olm\ra=0\big\}=\big\{u\in{\cal S}_+^2\;\big|\;u\olm=0\big\}=
\big\{\diag(0,a)\;\big|\;a\ge 0\big\}.
\end{equation*}
It follows from $\nabla\Phi(\ox)=\Big(\dfrac{\sub\Phi(\ox)}{\sub x_1},\dfrac{\sub\Phi(\ox)}{\sub x_2}\Big)=\big(\diag(1,0),\diag(0,1)\big)$ that
\begin{eqnarray*}
\ker\nabla\Phi(\ox)^*&=&\Big\{a=\begin{pmatrix}
a_{11}&a_{12}\\
a_{12}&a_{22}\\
\end{pmatrix}\in{\cal S}^2\;\Big|\;\Big(\Big\la a,\frac{\sub\Phi(\ox)}{\sub x_1}\Big\ra,\Big\la a,\frac{\sub\Phi(\ox)}{\sub x_2}\Big\ra\Big)=\nabla\Phi(\ox)^*a=0\Big\}\\
&=&
\Big\{a=\begin{pmatrix}
0&a_{12}\\
a_{12}&0\\
\end{pmatrix}\in{\cal S}^2\;\Big|\;a_{12}\in\R\Big\}.
\end{eqnarray*}
In this way we arrive at the representation
\begin{eqnarray}\label{db01}
K_{{\cal S}_+^2}(\oz,\olm)^*\cap\ker\nabla\Phi(\ox)^*&=&\Big\{\begin{pmatrix}
b_{11}&b_{12}\\
b_{12}&b_{22}\\
\end{pmatrix}\in{\cal S}^2\;\Big|\;b_{22}\le 0\Big\}\cap\Big\{\begin{pmatrix}
0&a_{12}\\
a_{12}&0\\
\end{pmatrix}\in{\cal S}^2\;\Big|\;a_{12}\in\R\Big\}\nonumber\\
&=&
\Big\{\begin{pmatrix}
0&a_{12}\\
a_{12}&0\\
\end{pmatrix}\in{\cal S}^2\;\Big|\;a_{12}\in\R\Big\},
\end{eqnarray}
which shows that the dual qualification condition \eqref{gf02} does not hold for SDP \eqref{sdp1}. On the other hand, we get from Example~\ref{ex1} that $\Lambda_{\ss{c}}(\ox)=\{\olm\}$. Let us now check that the multiplier mapping $M_{\ox}$ is not calm at $((0,0),\olm)$. Observe that $M_{\ox}$ admits the representation
\begin{eqnarray*}
M_{\ox}(v,w)&=&\Big\{\lm=\begin{pmatrix}
\lm_{11}&\lm_{12}\\
\lm_{12}&\lm_{22}\\
\end{pmatrix}\in{\cal S}^2_-\;\Big|\;v=\nabla_x L(\ox,\lm),\;\lm\in N_{{\cal S}_+^2}(w)\Big\}\\
&=&\Big\{\lm=\begin{pmatrix}
\lm_{11}&\lm_{12}\\
\lm_{12}&\lm_{22}\\
\end{pmatrix}\in{\cal S}^2_-\Big|\;v=(1+\lm_{11},\lm_{22}),\;\la\lm, w\ra=0\Big\}
\end{eqnarray*}
with $(v,w)\in\R^2\times{\cal S}^2$. Pick an arbitrary $t>0$ and define $v_t:=(-\frac{t^2}{2},-\frac{t^2}{2})$, $w_t:=\diag(0,0)$, and
$\lm_t:=\begin{pmatrix}
-1-\frac{t^2}{2}&\frac{t}{2}\\
\frac{t}{2}&-\frac{t^2}{2}\\
\end{pmatrix}$.
It is easy to see that $\lm_t\in M_{\ox}(v_t,w_t)\cap\B_t(\olm)$ when $t$ is sufficiently small. However, we have the limit calculation
$$
\lim_{t\dn 0}\frac{\|\lm_t-\olm\|}{\|v_t\|+\|w_t\|}=\lim_{t\dn 0}\frac{\sqrt{\frac{t^4}{2}+\frac{t^2}{2}}}{\frac{t^2}{\sqrt{2}}}=\infty,
$$
which shows that the mapping $M_{\ox}$ is not calm at $((0,0),\olm)$.}
\end{Example}\vspace*{-0.05in}

Observe to this end that in the NLP polyhedral framework we do not have the situation of Example~\ref{ex2}, since the calmness of $M_{\ox}$ is a direct consequence of the Hoffman lemma. In Section~\ref{sect5} we reveal a similar phenomenon telling us that $M_{\ox}$ is automatically calm in general nonpolyhedral systems under the strict complementarity condition formulated therein.\vspace*{-0.07in}

\begin{Remark}[\bf another characterization of uniqueness of Lagrange multipliers]\label{shapiro} {\rm In the case of optimization problems with the constraints $\Phi(x)\in\Th$ generated by convex cones $\Th$, Shapiro \cite[Proposition~2.1]{sh97} obtained a characterization of the uniqueness of Lagrange multipliers in the form
\begin{equation}\label{sh97}
\big(N_\Th(\Phi(\ox))+\R\olm\big)\cap\ker\nabla\Phi(\ox)^*=\{0\}.
\end{equation}
His result can be extended to the case of regular sets $\Th$ in the framework of Theorem~\ref{unique} by the following arguments. Assuming that the multiplier $\olm\in\Lambda(\ox)$ is unique, pick $\eta$ from the left-hand side of \eqref{sh97} and get $\eta=\lm+t\olm$ for some $\lm\in N_\Th(\Phi(\ox))$ and $t\in\R$. It follows from the regularity of $\Th$ that $\olm+\eta\in\Lambda(\ox)$ if $t\ge 0$ and that $\olm-\frac{1}{2t}\eta\in \Lambda(\ox)$ otherwise. This clearly contradicts the uniqueness of $\olm$. The converse implication can be also justified while arguing by contradiction. We see in the next section that the dual qualification condition \eqref{gf02} and the entire Theorem~\ref{unique} are very instrumental to derive complete characterizations of noncritical multipliers for \eqref{VS}. It seems not to be the case for condition \eqref{sh97}.}
\end{Remark}\vspace*{-0.3in}

\section{Characterizations of Noncritical Multipliers}
\sce\label{sect5}\vspace*{-0.1in}

In this section we establish the main result of the paper that gives us a complete characterization of noncriticality of Lagrange multipliers in general variational systems \eqref{VS}. Our previous result in  this direction \cite[Theorem~4.1]{ms17} addresses KKT systems of type \eqref{VS}
with $N_\Th$ replaced by the subdifferential mapping of a convex piecewise linear function. The proof therein is strongly based on the polyhedral structure of the latter systems and cannot be extended to a nonpolyhedral case. Here we develop a new approach that works for the general ${\cal C}^2$-cone reducible sets $\Th$.\vspace*{0.02in}

First we present several lemmas of their own interest.\vspace*{-0.1in}

\begin{Lemma}[\bf closed images under surjectivity] \label{axi0} Let $h\colon\Y\to\E$ be ${\cal C}^2$-smooth around $\oz$, and let $\nabla h(\oz)$ have full rank. Then $D\subset\E$ is closed if and only if $\nabla h(\oz)^*D$ has this property.
\end{Lemma}\vspace*{-0.17in}
\begin{proof} The `if' part comes as a direct consequence of the surjectivity condition $\ker\nabla h(\oz)^*=\{0\}$. The `only if' part follows from \cite[Lemma~1.18]{m06}.
\end{proof}\vspace*{-0.17in}

\begin{Lemma}[\bf propagation of closedness]\label{axi1}
Let the pair $(\ox,\olm)$ be a solution to the variational system \eqref{VS}, and let $\Th$ be ${\cal C}^2$-cone reducible at $\oz=\Phi(\ox)$ to a closed convex cone $C$. Then the following assertions are equivalent:\vspace{-0.2cm}
\begin{itemize}[noitemsep]
\item[\bf{(i)}] The set $K_\Th(\oz,\olm)^*-\big[K_\Th(\oz,\olm)^*\cap\ker\nabla\Phi(\ox)^*\big]$ is closed.
\item[\bf{(ii)}] The set $K_C(h(\oz),\bar\mu)^*-\big[K_C\big(h(\oz),\bar\mu\big)^*\cap\ker\nabla(h\circ\Phi)(\ox)^*\big]$ is
closed, where $h$ is taken from \eqref{red}, where $\bar\mu$ is a unique solution to \eqref{redl}, and where $K_C(h(\oz),\bar\mu):=T_C(h(\oz))\cap\{\bar\mu\}^\perp$ is the critical cone to $C$ at $h(\oz)$ for $\bar\mu$.
\end{itemize}
\end{Lemma}\vspace*{-0.17in}
\begin{proof} It follows from \eqref{gdf} that $DN_\Th(\oz,\olm)(0)=K_\Th(\oz,\olm)^*$. Thus the set in (i) can be equivalently represented as
\begin{equation*}
DN_\Th(\oz,\olm)(0)-\big[DN_\Th(\oz,\olm)(0)\cap\ker\nabla\Phi(\ox)^*\big].
\end{equation*}
Since $C$ is a closed convex cone with $h(\oz)=0\in C$, we conclude that $C$ is ${\cal C}^2$-cone reducible at $h(\oz)$ to itself in the sense of \eqref{red} with $h=I\colon\E\to\E$ being the identity mapping. This yields
\begin{equation}\label{io01}
DN_C\big(h(\oz),\bar\mu\big)(v)=N_{\ss{K_C(h(\oz),\bar\mu)}}(v)\;\mbox{ for all }\;v\in\E.
\end{equation}
Using the equivalent local representation \eqref{red} for $\Th$ and the surjectivity/full rank of $\nabla h(\oz)$, we deduce from \eqref{io01} and the second-order chain rule in \cite[Theorem~2]{go17} that
\begin{equation}\label{lm02}
DN_\Th(\oz,\olm)(u)= \nabla^2\la\bar\mu,h\ra(\oz)u+\nabla h(\oz)^*N_{\ss{K_C(h(\oz),\bar\mu)}}\big(\nabla h(\oz)u\big)\;\mbox{ for all }\;u\in\Y,
\end{equation}
which in turn implies the equalities
$$
DN_\Th(\oz,\olm)(0)=\nabla h(\oz)^*N_{\ss{K_C(h(\oz),\bar\mu)}}(0)=\nabla h(\oz)^*DN_C\big(h(\oz),\bar\mu\big)(0).
$$
The latter leads us to the representation
\begin{equation}\label{lm03}
\begin{array}{lll}
DN_\Th(\oz,\olm)(0)-\big[\big(DN_\Th(\oz,\olm)(0)\cap\ker\nabla\Phi(\ox)^*\big]\\=\nabla h(\oz)^*\Big\{DN_C\big(h(\oz),\bar\mu\big)(0)-\big[DN_C\big(h(\oz),\bar\mu\big)(0)\cap\ker\nabla(h\circ\Phi)(\ox)^*\big]\Big\}.
\end{array}
\end{equation}
Thus the claimed result amounts to saying that the following assertions are equivalent:\vspace{-0.2cm}
\begin{itemize}[noitemsep]
\item[\bf{(a)}] The set $DN_\Th(\oz,\olm)(0)-\big[\big(DN_\Th(\oz,\olm)(0)\cap\ker\nabla\Phi(\ox)^*\big]$ is closed.
\item[\bf{(b)}] The set $DN_C\big(h(\oz),\bar\mu\big)(0)-\big[DN_C\big(h(\oz),\bar\mu\big)(0)\cap\ker\nabla(h\circ\Phi)(\ox)^*\big]$
is closed.
\end{itemize}\vspace{-0.2cm}
Employing now \eqref{lm03} together with Lemma~\ref{axi0} readily verifies the equivalence between (a) and (b), and consequently between (i) and (ii).
\end{proof}\vspace*{-0.05in}

Consider next the set-valued mapping $S\colon\X\times\Y\tto\X\times\Y$ given by
\begin{equation}\label{mapS}
S(v,w):=\big\{(x,\lm)\in\X\times\Y\;\big|\;(v,w)\in G(x,\lm)\big\}\;\mbox{ for }\;(v,w)\in\X\times\Y,
\end{equation}
where the mapping $G$ is taken from \eqref{GKKT}. We can see that \eqref{mapS} defines the {\em solution map} to the {\em canonical perturbation} of the original variational system \eqref{VS}. The counterpart of \eqref{mapS} for the reduced generalized equation \eqref{rVS} is
\begin{equation}\label{rmapS}
S^r(v,w):=\big\{(x,\mu)\in\X\times\E\;\big|\;(v,w)\in G^r(x,\mu)\big\}\;\mbox{ with }\;(v,w)\in\X\times\E,
\end{equation}
where the corresponding mapping $G^r$ for \eqref{rVS} is defined by
\begin{equation}\label{RGKKT}
G^r(x,\mu):=\left[\begin{array}{c}
\Psi^r(x,\mu)\\
-(h\circ\Phi)(x)
\end{array}
\right]+\left[\begin{array}{c}
0\\
N_C^{-1}(\mu)
\end{array}
\right].
\end{equation}
The following lemma establishes the equivalence between an important stability property for the mappings $S$ and $S^r$ we introduced in \cite{ms17} under the name of {\em semi-isolated calmness}.\vspace*{-0.05in}

\begin{Lemma}[\bf propagation of semi-isolated calmness for solution mappings]\label{axi4} Let $(\ox,\olm)$ be a solution to the variational system \eqref{VS}, where $\Th$ is ${\cal C}^2$-cone reducible at $\oz=\Phi(\ox)$ to a closed convex cone $C$. Then the following assertions are equivalent: \vspace{-0.2cm}
\begin{itemize}[noitemsep]
\item[\bf{(i)}] There are numbers $\ve>0$ and $\ell\ge 0$ as well as neighborhoods $V$ of $0\in\X$ and $W$ of $0\in\Y$ such that for any
$(v,w)\in V\times W$ we have
\begin{equation}\label{upper}
S(v,w)\cap\B_\ve(\ox,\olm)\subset\big[\{\ox\}\times\Lm(\ox)\big]+\ell\big(\|v\|+\|w\|\big)\B.
\end{equation}
\item[\bf{(ii)}] There are numbers $\ve'>0$ and $\ell'\ge 0$ as well as neighborhoods $V$ of $0\in\X$ and $W$ of $0\in\E$ such that for any
$(v,w)\in V\times W$ we have
\begin{equation}\label{rupper}
S^r(v,w)\cap\B_{\ve'}(\ox,\omu)\subset\big[\{\ox\}\times\Lm^r(\ox)\big]+\ell'\big(\|v\|+\|w\|\big)\B.
\end{equation}
\end{itemize}
\end{Lemma}\vspace*{-0.17in}
\begin{proof}
Since $\nabla h(\oz)$ is surjective, there is a $\dd>0$ such that for any $z\in\B_\dd(\oz)$ the derivative $\nabla h(z)$ is surjective.
Pick $z\in U$ and find by \cite[Lemma~1.18]{m06} a constant $\kappa_z>0$ for which
$$
\kappa_z\|y\|\le\|\nabla h(z)^*y\|\;\mbox{ whenever }\;y\in\E.
$$
Denote $\bar\kappa:=\inf\{\kappa_z\;|\;z\in\B_{\dd/2}(\oz)\}$ and observe that $\bar\kappa>0$. Let us show then that
\begin{equation}\label{lm06}
\bar\kappa\|y\|\le\|\nabla h(z)^*y\|\;\mbox{ for all }\;z\in \B_{\dd/2}(\oz)\;\mbox{ and }\;y\in\E.
\end{equation}
Indeed, it follows from \cite[Lemma~1.18]{m06} that $\kappa_z=\inf\{\|\nabla h(z)^*y\|\;|\;\|y\|=1\}$ whenever $z\in\B_{\dd/2}(\oz)$. If $\bar\kappa=0$, we find a sequence of $z_k\in\B_{\dd/2}(\oz)$ with $\kappa_{z_k}\to 0$ as $k\to\infty$. This implies that there is a sequence of $y_k$ with $\|y_k\|=1$ such that
$$
\|\nabla h(z_k)^*y_k\|\le\kappa_{z_k}+k^{-1},\quad k\in\N.
$$
Passing to subsequences if necessary, assume without loss of generality that $z_k\to\tilde{z}$ and $y_k\to\tilde{y}$ with
$\tilde{z}\in\B_{\dd/2}(\oz)$ and $\|\tilde{y}\|=1$. Thus we arrive at $\nabla h(\tilde{z})^*\tilde{y}=0$, and hence $\tilde{y}=0$
due to the surjectivity of $\nabla h(\tilde{z})$. The obtained contradiction verifies \eqref{lm06}.

Assume now that (i) holds. Taken $\ve$ from (i), suppose without loss of generality that $\ell>0$ is a Lipschitz constant for the mappings
$\nabla h$ on $\B_{\ve}(\oz)$ and $\Phi$ on $\B_{\ve}(\ox)$. Let $M>0$ be an upper bound for the values of $\|\nabla h(\cdot)\|$ on $\B_{\ve}(\oz)$ and of $\|\nabla\Phi(\cdot)\|$ on $\B_{\ve}(\ox)$. It follows from \cite[Theorem~1.57]{m06} and the surjectivity of $\nabla h(\oz)$ that $h$ is metrically regular around $(\oz,0)$, i.e., there exist constants $\al>0$ and $\rho\ge 0$ such we have the estimate
\begin{equation}\label{lm08}
d\big(z;h^{-1}(y)\big)\le\rho\,\|h(z)-y\|\;\;\mbox{for all}\;\;(z,y)\in\B_\al(\oz)\times\B_\al(0).
\end{equation}
We can always suppose that $\B_\al(\oz)\subset{\cal O}$ with ${\cal O}$ taken from \eqref{red}. To prove the semi-isolated calmness of the mapping $S^r$ at $((0,0),(\ox,\omu))$, we claim that inclusion \eqref{rupper} holds with
\begin{equation}\label{lm10}
0<\ve'\le\min\Big\{\dfrac{\ve}{4\rho},\dfrac{\ve}{4\rho\ell\|\omu\|},\dfrac{\ve}{4\ell},\dfrac{\ve}{4\ell^2\|\omu\|},\dfrac{\al}{1+\ell^2}, \dfrac{\ve}{4M},\dfrac{\ve}{4},\dfrac{\al}{\ell},\frac{\dd}{2\rho+2\ell}\Big\},
\end{equation}
$V:=\B_{\ve'}(0)$, and $W:=\B_{\ve'}(0)$. To proceed, pick $(v,w)\in\B_{\ve'}(0)\times\B_{\ve'}(0)$ and $(x,\mu)\in S^r(v,w)\cap\B_{\ve'}(\ox,\omu)$ and then get the relationships
\begin{equation}\label{lm07}
v=\Psi^r(x,\mu)\;\mbox{ and }\;w+(h\circ\Phi)(x)\in N_C^{-1}(\mu).
\end{equation}
Let $y_w:=w+(h\circ\Phi)(x)$ and observe from \eqref{lm10} that $(\Phi(x),y_w)\in\B_\al(\oz)\times\B_\al(0)$.
Setting $z:=\Phi(x)$ and $y:=y_w$ in \eqref{lm08} gives us $z_w\in\Y$ such that
\begin{equation}\label{lm09}
\|\Phi(x)-z_w\|\le\rho\|w\|\;\mbox{ and }\;h(z_w)=y_w.
\end{equation}
This together with \eqref{rch} and \eqref{lm07} tells us that
$$
\begin{array}{lll}
v'=\Psi(x,\lm),&w'+\Phi(x)\in N_\Th^{-1}(\lm)\;\;\mbox{with}\\
\lm:=\nabla h(z_w)^*\mu,&w':=z_w-\Phi(x),\;\;v':=v+\nabla\Phi(x)^*(\nabla h(z_w)-\nabla h\big(\Phi(x))\big)^*\mu.
\end{array}
$$
Using \eqref{lm10}, we have the estimates
$$
\|z_w-\oz\|\le\|z_w-\Phi(x)\|+\|\Phi(x)-\Phi(\ox)\|\le\rho\|w\|+\ell\|x-\ox\|\le\min\Big\{\dfrac{\ve}{2},\dfrac{\ve}{2\ell\|\omu\|}\Big\},
$$
which yield in turn the following inequalities:
$$
\begin{array}{lll}
\|\lm-\olm\|&\le&\|\nabla h(z_w)\|\cdot\|\mu-\omu\|+\|\nabla h(z_w)-\nabla h(\oz)\|\cdot\|\omu\|\\
&\le&M\|\mu-\mu\|+\ell\|\omu\|\cdot\|z_w-\oz\|\le \dfrac{3\ve}{4}.
\end{array}
$$
This implies that $(x,\lm)\in S(v',w')\cap\B_{\ve}(\ox,\olm)$. It follows from (i) that there is a multiplier $\lm'\in\Lambda(\ox)$ such that $\|x-\ox\|+\|\lm-\lm'\|\le\ell(\|v'\|+\|w'\|)$. Using \eqref{lare} gives us $\mu'\in\Lambda^r(\ox)$ such that $\lm'=\nabla h(\oz)^*\mu'$. Then we get from \eqref{lm10} that $z_w\in\B_{\dd/r}(\oz)$, which ensures by \eqref{lm06} that
\begin{eqnarray*}
\bar\kappa\|\mu-\mu'\|&\le&\|\nabla h(z_w)^*\mu-\nabla h(z_w)^*\mu'\|\\
&\le&\|\nabla h(z_w)^*\mu-\nabla h(\oz)^*\mu'\|+\|\nabla h(z_w)-\nabla h(\oz)\|\cdot\|\mu\|\\
&\le&\|\lm-\lm'\|+\ell\|z_w-\oz\|(\ve+\|\omu\|).
\end{eqnarray*}
This allows us to obtain the relationships
$$
\begin{array}{lll}
\|x-\ox\|+\|\mu-\omu\|&\le\|x-\ox\|+\dfrac{1}{\bar\kappa}\|\lm -\lm'\|+\dfrac{\ell(\ve+\|\omu\|)}{\bar\kappa}\|z_w-\oz\|\\
&\le\|x-\ox\|+\dfrac{1}{\bar\kappa}\|\lm-\lm'\|+\dfrac{\ell(\ve+\|\omu\|)}{\bar\kappa}\big(\rho\|w\|+\ell\|x-\ox\|\big)\\
&\le\max\Big\{\dfrac{1}{\bar\kappa},1+\dfrac{\ell^2(\ve+\|\omu\|)}{\bar\kappa}\Big\}\big(\|x-\ox\|+\|\lm-\lm'\|\big)+\dfrac{\ell\rho (\ve+\|\omu\|)}{\bar\kappa}\|w\|\\
&\le\max\Big\{\dfrac{1}{\bar\kappa},1+\dfrac{\ell^2(\ve+\|\omu\|)}{\bar\kappa}\Big\}\ell\,\big(\|v'\|+\|w'\|\big)+\dfrac{\ell\rho(e+\|\omu\|)}
{\bar\kappa}\|w\|\\
&\le\max\Big\{\dfrac{1}{\bar\kappa},1+\dfrac{\ell^2(\ve+\|\omu\|)}{\bar\kappa}\Big\}\ell\,\big(\|v\|+M(\ve+\|\omu\|)\ell\rho\|w\|+\rho\|w\|\big)\\
&\;\;\;+\dfrac{\ell\rho(\ve+\|\omu\|)}{\bar\kappa}\|w\|,
\end{array}
$$
which therefore verify the claimed inclusion \eqref{lm05}.

Suppose next that the mapping $S^r$ is semi-isolatedly calm at $((0,0),(\ox,\bar\mu))$ and thus find constants $\ell'\ge 0$ and $\ve'>0$ for which \eqref{rupper} is satisfied. We can always assume that $\ell$ is a Lipschitz constant for the mappings $\nabla h$ on $\B_{\ve'}(\oz)$ and $\Phi$ on $\B_{\ve'}(\ox)$ and that $M$ is an upper bound for $\|\nabla\Phi(\cdot)\|$ on $\B_{\ve'}(\ox)$. To prove \eqref{upper}, take $\ve>0$ such that
\begin{equation}\label{lm11}
\ve\le\min\Big\{\dfrac{\ve'}{4(\ell+1)},\dfrac{\dd}{4(\ell+1)},\dfrac{\bar\kappa\ve'}{2(1+(\ell+\ell^2)\|\omu\|)},\dfrac{\ve'}{4}\Big\},
\end{equation}
where $\dd$ is taken from \eqref{lm06}, and suppose that $\B_{\ve'}(\oz)\subset{\cal O}$ with ${\cal O}$ taken from \eqref{red}. Picking $(v,w)\in\B_{\ve}(0)\times\B_{\ve}(0)$, we get $(x,\lm)\in S(v,w)\cap\B_{\ve}(\ox,\olm)$ and hence
$$
v=\Psi(x,\lm)\;\mbox{ and }\;w+\Phi(x)\in N_\Th^{-1}(\lm).
$$
Let $z_w:=w+\Phi(x)$ and deduce from \eqref{lm11} that $z_w\in\B_{\ve'}(\oz)\subset{\cal O}$. This tells us by \eqref{rch} that
$$
\lm:=\nabla h(z_w)^*\mu\;\mbox{ for some }\;\mu\in N_C\big(h(z_w)\big),
$$
which ensures therefore that
$$
\begin{array}{lll}
v'=\Psi^r(x,\mu),\;\; w'+ h\big(\Phi(x)\big)\in N_C^{-1}(\mu)\;\;\mbox{with}\\
w':=\nabla h\big(\Phi(x)\big)w+o(\|w\|),\;\;v':=v+\nabla\Phi(x)^*\big(\nabla h(z_w)-\nabla h(\Phi(x))\big)^*\mu.
\end{array}
$$
It follows from \eqref{lm11} that $z_w\in\B_{\dd/2}(\oz)$, and thus \eqref{lm06} leads us to the estimates
$$
\begin{array}{lll}
\|\mu-\omu\|&\le&\dfrac{1}{\bar\kappa}\|\nabla h(z_w)^*(\mu-\omu)\|\\
&\le&\dfrac{1}{\bar\kappa}\big(\|\nabla h(z_w)^*\mu-\nabla h(\oz)^*\omu\|+\|\nabla h(z_w)-\nabla h(\oz)\|\cdot\|\omu\|\big)\\
&\le&\dfrac{1}{\bar\kappa}\big(\|\lm-\olm\|+\ell\|\omu\|\cdot\|z_w-\oz\|\big)\\
&\le&\dfrac{1}{\bar\kappa}\big(\|\lm-\olm\|+\ell\|\omu\|(\|w\|+\ell\|x-\ox\|)\big)\le\dfrac{\ve(1+(\ell+\ell^2)\|\omu\|)}{\bar\kappa}
\le\dfrac{\ve'}{2},
\end{array}
$$
which yield $(x,\mu)\in S^r(w',v')\cap\B_{\ve'}(\ox,\omu)$. Appealing now to \eqref{rupper} gives us $\mu'\in\Lambda^r(\ox)$ such that $\|x-\ox\|+\|\mu-\mu'\|\le \ell'(\|v'\|+\|w'\|)$. By \eqref{lare} we find $\lm'\in\Lambda (\ox)$ with $\lm'=\nabla h(\oz)^*\mu'$ and
\begin{eqnarray*}
\|\lm-\lm'\|&\le&\|\nabla h(z_w)-\nabla h(\oz)\|\cdot\|\mu\|+\|\nabla h(\oz)\|\cdot\|\mu-\mu'\|\\
&\le&(\|\omu\|+\ve')\ell\|z_w-\oz\|+\|\nabla h(\oz)\|\cdot\|\mu-\mu'\|\\
&\le&(\|\omu\|+\ve')\ell(\|w\|+\ell\|x-\ox\|)+\|\nabla h(\oz)\|\cdot\|\mu-\mu'\|.
\end{eqnarray*}
Therefore we arrive at the inequalities
\begin{eqnarray*}
\|x-\ox\|+\|\lm-\lm'\|&\le&\|x-\ox\|+(\|\omu\|+\ve')\ell(\|w\|+\ell\|x-\ox\|)+ \|\nabla h(\oz)\|\cdot\|\mu-\mu'\|\\
&\le&\max\big\{1+\ell^2(\|\omu\|+\ve'),\|\nabla h(\oz)\|\big\}\big(\|x-\ox\|+\|\mu-\mu'\|\big)+(\|\omu\|+\ve')\ell\|w\|\\
&\le&\max\big\{1+\ell^2(\|\omu\|+\ve'),\|\nabla h(\oz)\|\big\}\ell'\big(\|v'\|+\|w'\|\big)+ (\|\omu\|+\ve')\ell\|w\|\\
&\le&\max\big\{1+\ell^2(\|\omu\|+\ve'),\|\nabla h(\oz)\|\big\}\ell'\big(M\|w\|+\|o(\|w\|)\|+\|v\|\\
&&+M\ell(\|\omu\|+\ve')\|w\|\big)+(\|\omu\|+\ve')\ell\|w\|,
\end{eqnarray*}
which verify \eqref{lm04} and thus complete the proof of this lemma.
\end{proof}\vspace*{-0.1in}

Next we establish relationships between the calmness property \eqref{calm-def} for the original system \eqref{VS} and its reduced counterpart \eqref{red}. To proceed, pick a stationary point $\ox$ from \eqref{stat} and define the {\em reduced multiplier mapping} $M^r_{\ox}\colon\X\times\E\tto\E$ by
\begin{equation}\label{rlagmap}
M^r_{\ox}(v,w):=\big\{\mu\in\E\;\big|\;(v,w)\in G^r(\ox,\mu)\big\}\;\mbox{ with }\;(v,w)\in\X\times\E.
\end{equation}\vspace*{-0.3in}

\begin{Lemma}[\bf propagation of calmness for multiplier mappings]\label{axi3} Let $(\ox,\olm)$ be a solution to the variational system \eqref{VS}, where $\Th$ is ${\cal C}^2$-cone reducible at $\oz=\Phi(\ox)$ to a closed convex cone $C$. Then the calmness of the mapping $M_{\ox}$ from \eqref{lagmap} at $((0,0),\olm)$ is equivalent to that of the mapping $M^r_{\ox}$ from \eqref{rlagmap} at $((0,0),\bar\mu)$, where $\bar\mu$ is a unique solution to \eqref{redl}.
\end{Lemma}\vspace*{-0.15in}
\begin{proof} The calmness property of $M_{\ox}$ at $((0,0),\olm)$ gives us $\ell\ge 0$ and $\ve>0$ with
\begin{equation}\label{lm04}
M_{\ox}(v,w)\cap\B_{\ve}(\olm)\subset M_{\ox}(0,0)+\ell\big(\|v\|+\|w\|\big)\B\;\mbox{ for all }\;(v,w)\in\B_{\ve}(0,0).
\end{equation}
To verify the calmness of $M^r_{\ox}$ at $((0,0),\bar\mu)$, we show that
\begin{equation}\label{lm05}
M^r_{\ox}(v,w)\cap\B_{\ve'}(\bar\mu)\subset M^r_{\ox}(0,0)+\ell'\big(\|v\|+\|w\|\big)\B\;\mbox{ whenever }\;(v,w)\in\B_{\ve'}(0,0)
\end{equation}
for $\ve':=\min\big\{\dfrac{\ve}{\|\nabla h(\oz)\|},\dfrac{\ve}{2}\big\}$ and $\ell':=\dfrac{\ell}{\bar\kappa}$ with $\bar\kappa$ taken from \eqref{lm06}. To proceed, pick $(v,w)\in\B_{\ve'}(0,0)$ and $(v,w)\in M^r_{\ox}(v,w)\cap\B_{\ve'}(\bar\mu)$ telling us that
$$
v=\Psi^r(\ox,\mu)\;\mbox{ and }\;w+h(\oz)\in N_C^{-1}(\mu).
$$
Since $h(\oz)=0$, we have $N_C(y)\subset C^*=N_C(h(\oz))$ for any $y\in\E$. Denoting $\lm:=\nabla h(\oz)^*\mu$, deduce from \eqref{rch} that the above conditions yield
$$
v=\Psi(\ox,\lm)\;\mbox{ and }\;\lm\in N_\Th(\oz),
$$
and thus $\lm\in M_{\ox}(v,0)$. It follows from $\mu\in\B_{\ve'}(\bar\mu)$ that $\lm\in\B_{\ve}(\olm)$. Combining this with \eqref{lm04}, we find $\lm'\in M_{\ox} (0,0)=\Lambda(\ox)$ such that $\|\lm-\lm'\|\le\ell\|v\|$. Invoking \eqref{lare} gives us $\mu'\in\Lambda^r(\ox)=M^r_{\ox}(0,0)$ with $\lm'=\nabla h(\oz)^*\mu'$. Remembering \eqref{lm06}, we arrive at the relationships
$$
\kappa\|\mu-\mu'\|\le\|\nabla h(\oz)^*\mu-\nabla h(\oz)^*\mu'\|=\|\lm-\lm'\|\le\ell\|v\|,
$$
which justify the claimed inclusion \eqref{lm05}.

Assume now that the mapping $M^r_{\ox} $ is calm at $((0,0),\bar\mu)$ and find constants $\ell'\ge 0$ and $\ve'>0$ for which \eqref{lm05} is satisfied. To prove \eqref{lm04} for the mapping $M_{\ox}$, select $\ve>0$ so that
$$
\ve\le\min\Big\{\dfrac{\ve'}{4},\dfrac{\ve'}{4\ell\|\nabla\Phi(\ox)\|(\|\bar\mu\|+\ve')},\dfrac{\bar\kappa\ve'}{4},\dfrac{\bar\kappa\ve'}{4\ell
(\|\bar\mu\|+\ve)}\Big\},
$$
where $\ell$ is a Lipschitz constant for $\nabla h$ around $\oz$. Picking $(v,w)\in\B_{\ve}(0,0)$ and $\lm\in M_{\ox}(v,w)\cap\B_{\ve}(\olm)$, we arrive at the conditions
$$
v=\Psi(\ox,\lm)\;\mbox{ and }\;w+\oz\in N_\Th^{-1}(\lm).
$$
Suppose without loss of generality that $w+\oz\in{\cal O}$, where the neighborhood ${\cal O}$ is taken from \eqref{red}. It allows us to deduce from \eqref{rch} that $\lm=\nabla h(w+\oz)^*\mu$ for some $\mu\in N_C(h(w+\oz))\subset N_C(h(\oz))$, and therefore to get
$$
v+\nabla\Phi(\ox)^*\big(\nabla h(\oz)-\nabla h(w+\oz)\big)^*\mu=\Psi^r(\ox,\mu)\;\mbox{ and }\;h(\oz)\in N_C^{-1}(\mu).
$$
This means that $\mu\in M^r_{\ox}(v',0)$ with $v'=v+\nabla\Phi(\ox)^*\big(\nabla h(\oz)-\nabla h(w+\oz)\big)^*\mu$. By using \eqref{lm06} and the selection of $\ve$ we obtain the inequalities
\begin{eqnarray*}
\|\mu-\bar\mu\|&\le&\dfrac{1}{\bar\kappa}\|\nabla h(\oz)^*(\mu-\omu)\|\\
&\le&\dfrac{1}{\bar\kappa}\|\lm-\olm\|+\dfrac{\|\nabla h(w+\oz)-\nabla h(\oz)\|(\|\bar\mu\|+\ve)}{\bar\kappa}\le\dfrac{\ve'}{4}+\dfrac{\ve'}{4}<\ve',
\end{eqnarray*}
which show that $\mu\in M^r_{\ox}(v',0)\cap\B_{\ve'}(\bar\mu)$ with $v'$ satisfying
$$
\|v'\|\le\|v\|+\ell\|\nabla\Phi(\ox)\|\cdot\|w\|\cdot\|\mu\|\le\dfrac{\ve'}{4}+\ell\|\nabla\Phi(\ox)\|\cdot\|w\|(\|\bar\mu\|+\ve')\le
\dfrac{\ve'}{4}+\dfrac{\ve'}{4}<\ve'.
$$
Appealing now to \eqref{lm05} gives us $\mu'\in M^r_{\ox}(0,0)=\Lambda^r(\ox)$ with $\|\mu -\mu'\|\le\ell'\|v'\|$. Employing \eqref{lare} again, we find $\lm'\in\Lambda(\ox)=M_{\ox}(0,0)$ such that $\lm'=\nabla h(\oz)^*\mu'$ and
\begin{eqnarray*}
\|\lm-\lm'\|&=&\|\nabla h(w+\oz)^*\mu-\nabla h(\oz)^*\mu'\|\\
&\le&\|\nabla h(w+\oz)-\nabla h(\oz)\|\cdot\|\mu\|+\|\nabla h(\oz)\|\cdot\|\mu-\mu'\|\\
&\le&(\|\bar\mu\|+\ve')\ell\|w\|+\ell'\|\nabla h(\oz)\|\cdot\|v'\|\\
&\le&(\|\bar\mu\|+\ve')\ell\|w\|+\ell'\|\nabla h(\oz)\|\big(\|v\|+\ell\|\nabla\Phi(\ox)\|(\|\bar\mu\|+\ve')\|w\|\big),
\end{eqnarray*}
which verifies \eqref{lm04} and thus completes the proof.
\end{proof}\vspace*{-0.1in}

The last lemma in this section establishes an equivalence between noncriticality of Lagrange multipliers of the original and reduced systems. \vspace*{-0.05in}

\begin{Lemma}[\bf propagation of noncriticality]\label{axi5} Let $(\ox,\olm)$ be a solution to the variational system \eqref{VS}, and let $\Th$ be ${\cal C}^2$-cone reducible at $\oz=\Phi(\ox)$ to the closed convex cone $C$. Then the Lagrange multiplier $\olm\in\Lambda(\ox)$ from \eqref{laset} is noncritical for \eqref{VS} if and only if the unique solution $\bar\mu\in\Lambda^r(\ox)$ to \eqref{redl} is noncritical for \eqref{rVS}.
\end{Lemma}\vspace*{-0.15in}
\begin{proof} Employing the classical chain rule, we get
\begin{eqnarray*}
\nabla^2\la\bar\mu,h\circ\Phi\ra(\ox)&=&\nabla\big(\nabla(h\circ\Phi)(\ox)^*\omu\big)=\nabla\big[\nabla\Phi(\cdot)^*\nabla h\big(\Phi(\cdot)\big)^*\omu\big]\Big|_{x=\ox}\\
&=&\nabla\big[\nabla\Phi(\cdot)^*\nabla h\big(\Phi(\ox)\big)^*\omu\big]\Big |_{x=\ox}+\nabla\big[\nabla\Phi(\ox)^*\nabla h\big(\Phi(\cdot)\big)^*\omu\big]\Big|_{x=\ox}\\&=&\nabla\big[\nabla\Phi(\cdot)^*\big(\nabla h(\oz)^*\bar\mu\big)\big]\Big|_{x=\ox}+\nabla\Phi(\ox)^*\nabla^2\la\bar\mu,h\ra(\oz)\nabla\Phi(\ox)\\
&=&\nabla^2\la\olm,\Phi\ra(\ox)+\nabla\Phi(\ox)^*\nabla^2\la\bar\mu,h\ra(\oz)\nabla\Phi(\ox).
\end{eqnarray*}
Combining this with \eqref{crc}, \eqref{lm02}, and \eqref{io01} yields the relationships
\begin{eqnarray*}
&&\nabla_{x}\Psi(\ox,\olm)\xi+\nabla\Phi(\ox)^*DN_\Th\big(\Phi(\ox),\olm\big)\big(\nabla\Phi(\ox)\xi\big)\\
&=&\nabla_{x}\Psi(\ox,\olm)\xi+\nabla\Phi(\ox)^*\Big\{\nabla^2\la\bar\mu,h\ra(\oz)\nabla\Phi(\ox)\xi+\nabla h(\oz)^* N_{\ss{K_C(h(\oz),\bar\mu)}}\big(\nabla h(\oz)\nabla\Phi(\ox)\xi\big)\Big\}\\
&=&\nabla f(\ox)\xi+\nabla^2\la\olm,\Phi\ra(\ox)\xi+\nabla\Phi(\ox)^*\nabla^2\la\bar\mu, h\ra(\oz)\nabla\Phi(\ox)\xi\\
&&+\nabla(h\circ\Phi)(\ox)^*DN_C(h(\oz),\bar\mu)\big(\nabla(h\circ\Phi)(\ox)\xi\big)\\
&=&\nabla f(\ox)\xi+\nabla^2\la\bar\mu,(h\circ\Phi)\ra(\ox)\xi+\nabla(h\circ\Phi)(\ox)^*DN_C(h(\oz),\bar\mu)\big(\nabla(h\circ\Phi)(\ox)\xi\big)\\
&=&\nabla_{x}\Psi^r(\ox,\bar\mu)\xi+\nabla(h\circ\Phi)(\ox)^*DN_C(h(\oz),\bar\mu)\big(\nabla(h\circ\Phi)(\ox)\xi\big),
\end{eqnarray*}
which justify the claimed equivalence for noncritical Lagrange multipliers.
\end{proof}\vspace{-0.2cm}

Now we are ready to establish the main result of the paper that provides a complete characterization of noncriticality
of Lagrange multipliers for nonpolyhedral variational systems \eqref{VS}.\vspace*{-0.05in}

\begin{Theorem}[\bf characterizations of noncritical Lagrange multipliers]\label{uplip} Let $(\ox,\olm)$ be a solution to the variational system \eqref{VS}. Consider the following properties of \eqref{VS} and the solution map $S$ taken from \eqref{mapS}:\vspace{-0.2 cm}
\begin{itemize}[noitemsep]
\item[\bf{(i)}] The Lagrange multiplier $\olm\in\Lambda(\ox)$ from \eqref{laset} is noncritical for \eqref{VS}.
\item[\bf{(ii)}] There are numbers $\ve>0$, $\ell\ge 0$ and neighborhoods $V$ of $0\in\X$ and $W$ of $0\in\Y$ such that for any
$(v,w)\in V\times W$ the semi-isolated calmness inclusion \eqref{upper} holds.
\item[\bf{(iii)}] There are numbers $\ve>0$ and $\ell\ge 0$ such that the estimate
\begin{equation}\label{subr}
\|x-\ox\|+d\big(\lm;\Lm(\ox)\big)\le\ell\big(\|\Psi(x,\lm)\|+d\big(\Phi(x);N^{-1}_{\Th}(\lm)\big)\big)
\end{equation}
is satisfied for all pairs $(x,\lm)\in\B_\ve(\ox,\olm)$.
\end{itemize}\vspace*{-0.05in}
Then we have the assertions:\vspace{-0.2cm}
\begin{itemize}[noitemsep]
\item[\bf{(a)}] Implications {\rm(iii)}$\Longleftrightarrow${\rm(ii)}$\Longrightarrow${\rm(i)} always fulfill.
\item[\bf{(b)}] If $\Th$ is ${\cal C}^2$-cone reducible at $\oz=\Phi(\ox)$ to a closed convex cone $C$, if the set
\begin{equation}\label{lm15}
K_\Th(\oz,\olm)^*-\big[K_\Th(\oz,\olm)^*\cap\ker\nabla\Phi(\ox)^*\big]
\end{equation}
is closed, and if the Lagrange multiplier mapping $M_{\ox}$ from \eqref{lagmap} is calm at $((0,0),\olm)$, then the converse implication {\rm(i)}$\Longrightarrow${\rm(ii)} is also satisfied.
\end{itemize}\vspace{-0.2cm}
\end{Theorem}\vspace*{-0.2in}
\begin{proof} The equivalence between (ii) and (iii) can be verified similarly to \cite[Theorem~4.1]{ms17}. To prove (ii)$\Longrightarrow$(i), it suffices to show that \eqref{cri2} holds. Pick $(\xi,\eta)\in\X\times\Y$ satisfying
$(0,0)\in DG\big((\ox,\olm),(0,0)\big)(\xi,\eta)$ and get $\big((\xi,\eta),(0,0)\big)\in T_{\ss{\gph G}}\big((\ox,\olm),(0,0)\big )$. By the definition of the graphical derivative, find sequences $t_k\downarrow 0$ and $\big((\xi_k,\eta_k),(v_{k},w_{k})\big)\to\big((\xi,\eta),(0,0)\big)$ with
$$
\big((\ox,\olm),(0,0)\big)+t_k\big((\xi_k,\eta_k),(v_k,w_k)\big)\in\gph G\;\mbox{ for all }\;k\in\N.
$$
Remembering the definition of $S$ in \eqref{mapS} gives us the inclusions
$$
(\ox+t_k\xi_k,\olm+t_k\eta_k)\in S(t_k v_{k},t_k w_{k}),\quad k\in\N.
$$
It follows from \eqref{upper} that for all $k$ sufficiently large we have
$$
t_k\|\xi_k\|=\|x_t-\ox\|\le\ell t_k\big(\|v_{t}\|+\|w_{t}\|\big).
$$
Divining there by $t_k$ and then letting $k\to\infty$ imply that $\xi=0$, and thus (a) holds.

Turning to (b), we appeal to Lemma~\ref{axi5}, which tells us that $\omu$ from \eqref{redl} is a noncritical multiplier for \eqref{rVS}. Let us show that the mapping $S^r$ from \eqref{rmapS} is semi-isolatedly calm at $\big((0,0),(\ox,\omu)\big)$, i.e.,
inclusion \eqref{rupper} holds for some constants $\ve'>0$ and $\ell'\ge 0$ and for some neighborhoods $V$ of $0\in\X$ and $W$ of $0\in\E$.
To furnish this, we first verify the following result.\\[1ex]
{\bf Claim:} {\em There are numbers $\ve'>0$, $\ell'\ge 0$ and neighborhoods $V$ of $0\in\X$ and $W$ of $0\in\E$
such that for any $(v,w)\in V\times W$ and any $(x_{vw},\mu_{vw})\in S^r(v,w)\cap\B_{\ve'}(\ox,\omu)$ we have the estimate}
\begin{equation}\label{upper2}
\|x_{vw}-\ox\|\le\ell'\big(\|v\|+\|w\|\big).
\end{equation}
To prove this claim, suppose on the contrary that \eqref{upper2} fails, i.e., for any $k\in\N$ there are
$(v_{k},w_{k})\in\B_{1/k}(0)\times\B_{1/k}(0)$ and $(x_k,\mu_k)\in S^r(v_{k},w_{k})\cap\B_{1/k}(\ox,\omu)$ satisfying
$$
\frac{\|x_k-\ox\|}{\|v_{k}\|+\|w_{k}\|}\to\infty\;\mbox{ as }\;k\to\infty,
$$
which yields $v_{k}=o(\|x_k-\ox\|)$ and $w_{k}=o(\|x_k-\ox\|)$. Letting $y_k:=(h\circ\Phi)(x_k)+w_{k}$, observe from (\ref{rmapS}) that $(y_k,\mu_k)\in\gph N_C$. We know from Lemma~\ref{axi3} that the calmness property for $M_{\ox}$ at $((0,0),\olm)$ amounts to that for $M^r_{\ox}$ at $((0,0),\omu)$. The latter is equivalent to the metric subregularity of $(M^r_{\ox})^{-1}$ at $(\omu,(0,0))$, which gives us  $\rho\ge 0$ and $\al>0$ such that
\begin{equation}\label{lm13}
d\big(\mu;\Lambda^r(\ox)\big)\le\rho\big(\|\Psi^r(\ox,\mu)\|+d\big(h(\oz);N_C^{-1}(\mu)\big)\big)\;\mbox{ for all }\;\mu\in\B_\al(\omu).
\end{equation}
This together with $h(\oz)=0$ allows us to get for all $k$ sufficiently large the estimates
\begin{eqnarray}\label{lm14}
d\big(\mu_k;\Lm^r(\ox)\big)&\le&\rho\big(\|\Psi^r(\ox,\mu_k)\|+d\big(h(\oz);N_C^{-1}(\mu_{k})\big)\nonumber\\
&=&\rho\big(\|f(\ox)+\nabla(h\circ\Phi)(\ox)^*\mu_{k}\|\big)\nonumber\\
&\le&\rho\big(\|\mu_{k}\|\cdot\|\nabla(h\circ\Phi)(x_{k})-\nabla(h\circ\Phi)(\ox)\|+\|\nabla(h\circ\Phi)(x_{k})^*\mu_{k}+f(x_{k})\|\nonumber\\
&&+\|f(x_{k})-f(\ox)\|\big)\nonumber\\
&\le&\rho\big(\ell''\|\mu_k\|\cdot\|x_{k}-\ox\|+\|v_{k}\|+\ell''\|x_{k}-\ox\|\big),
\end{eqnarray}
where $\ell''$ is a calmness constant for the mappings $f$ and $\nabla(h\circ\Phi)$ at $\ox$. Thus there is $\mu'_k\in\Lambda^r(\ox)$ such that the sequence $\dfrac{\mu_k-\mu'_k}{\|x_k-\ox\|}$ is bounded and so contains a convergent subsequence
\begin{equation}\label{seq1}
\eta_k:=\frac{\mu_k-\mu'_k}{\|x_k-\ox\|}\to\tet\;\mbox{ as }\;k\to\infty\;\mbox{ with some }\;\tet\in\E.
\end{equation}
Passing to a subsequence if necessary, we get that
\begin{equation}\label{xi}
\xi_k:=\frac{x_k-\ox}{\|x_k-\ox\|}\to\xi\;\mbox{ as }\;k\to\infty\;\mbox{ with some }\;0\ne\xi\in\X.
\end{equation}
Denote $t_k:=\|x_k -\ox\|$ and deduce from $(x_k,\mu_k)\in S^r(v_{k},w_{k})$ that
$$
o(t_k)=\Psi^r(x_{k},\mu_{k})\;\mbox{ and }\;\mu_k\in N_C(y_k).
$$
Taking this into account and using \eqref{seq1} lead us to
\begin{eqnarray*}
o(t_k)=v_{k}&=&\Psi^r(x_{k},\mu_{k})=\Psi^r(x_{k},\omu)-\Psi(\ox,\omu)+\nabla(h\circ \Phi)(x_k)^*(\mu_k-\omu)\\
&=&\nabla_{x}\Psi^r(\ox,\omu)(x_k-\ox)+\nabla(h\circ \Phi)(\ox)^*(\mu_k-\omu)+o(t_k)\\
&=&\nabla_{x}\Psi^r(\ox,\omu)(x_k-\ox)+\nabla(h\circ \Phi)(\ox)^*(\mu_k-\mu'_k)+o(t_k),
\end{eqnarray*}
which in turn yields the equality
\begin{equation}\label{up10}
\nabla_{x}\Psi^r(\ox,\omu)\xi+\nabla(h\circ\Phi)(\ox)^*\tet=0.
\end{equation}
Since $C$ is a closed convex cone, it follows from $(y_k,\mu_k)\in\gph N_C$ that $y_k\in C$ and $\la y_k,\mu_k\ra=0$. The latter together with $h(\oz)=0$ leads us to
$$
0=\la y_k,\mu_k\ra=\big\la w_k+h\big(\Phi(x_k)\big)-h(\oz),\mu_k\big\ra=\big\la\nabla(h\circ\Phi)(\ox)(x_k-\ox)+o(t_k),\mu_k\big\ra,
$$
and hence $\la\nabla(h\circ\Phi)(\ox)\xi,\omu\ra=0$. We have that $(h\circ\Phi)(\ox)+t_k\big[\nabla(h\circ\Phi)(\ox)\xi_k+o(t_k)/t_k\big]\in C$, which implies that$\nabla(h\circ\Phi)(\ox)\xi\in T_C(h(\oz))$ and so
\begin{equation}\label{tr03}
\nabla(h\circ\Phi)(\ox)\xi\in T_C\big(h(\oz)\big)\cap\{\omu\}^\bot\cap\rge\nabla(h\circ\Phi)(\ox)=K_C\big(h(\oz),\bar\mu\big)\cap\rge\nabla(h\circ\Phi)(\ox).
\end{equation}
It follows from $\mu'_k,\omu\in\Lambda^r(\ox)$ and $\mu_k\in N_C(y_k)$ that
\begin{eqnarray*}
\mu_k-\mu'_k&=&\mu_k-\omu+\omu-\mu'_k\in C^*+\R\omu-\big[\ker\nabla(h\circ\Phi)(\ox)^*\cap(C^*+\R\omu)\big]\\
&=&N_C\big(h(\oz)\big)+\R\omu-\big[\ker\nabla(h\circ\Phi)(\ox)^*\cap\big(N_C(h(\oz))+\R\omu\big)\big]\\
&\subset&\cl\big(N_C(h(\oz))+\R\omu\big)-\big[\ker\nabla(h\circ\Phi)(\ox)^*\cap\cl\big(N_C(h(\oz))+\R\omu\big)\big]\\
&=&K_C\big(h(\oz),\bar\mu\big)^*-\big[\ker\nabla(h\circ\Phi)(\ox)^*\cap\big(K_C\big(h(\oz),\bar\mu\big)\big)^*\big]\\
&=&\big(K_C\big(h(\oz),\bar\mu\big)\cap{\cal D}^*\big)^*,
\end{eqnarray*}
where ${\cal D}:=-\big[\ker\nabla(h\circ\Phi)(\ox)^*\cap\big(K_C(h(\oz),\bar\mu)\big)^*\big]$, and where the last equality comes from the closedness assumptions \eqref{lm15}, Lemma~\ref{axi1}, and \cite[Proposition~20]{bbw}. This leads us to
\begin{equation}\label{mn01}
\tilde{\eta}\in\big(K_C\big(h(\oz),\bar\mu\big)\cap{\cal D}^*\big)^*.
\end{equation}
On the other hand, we have $\rge\nabla(h\circ\Phi)(\ox)\subset{\cal D}^*$, which together with \eqref{tr03} yields
\begin{equation}\label{tr04}
\nabla(h\circ\Phi)(\ox)\xi\in K_C\big(h(\oz),\bar\mu\big)\cap{\cal D}^*.
\end{equation}
Remember that $\mu'_k\in N_C(h(\oz))$ and $\mu_k\in N_C(y_k)$. It follows from the monotonicity of normal cone mappings to convex sets that
$$
0\le\Big\la\frac{\mu_k-\mu'_k}{t_k},\frac{y_k-h(\oz)}{t_k}\Big\ra.
$$
This implies therefore that
\begin{equation}\label{mn02}
\la\nabla(h\circ \Phi)(\ox)\xi\tilde{\eta}\ra\ge 0.
\end{equation}
Taking this into account together with \eqref{mn01} and \eqref{tr04} implies that
\begin{equation}\label{mn03}
\tilde\eta\in N_{\ss{K_C(h(\oz),\bar\mu)\cap{\cal D}^*}}\big(\nabla(h\circ\Phi)(\ox)\xi\big).
\end{equation}
Appealing again the intersection rule from \cite[Proposition~20]{bbw} to \eqref{mn03} gives us
$$
\tilde\eta\in N_{\ss{K_C(h(\oz),\bar\mu)}}\big(\nabla(h\circ\Phi)(\ox)\xi\big)+N_{\ss{{\cal D}^*}}\big(\nabla(h\circ\Phi)(\ox)\xi\big).
$$
Thus there exist vectors $\eta\in N_{\ss{K_C(h(\oz),\bar\mu)}}\big(\nabla(h\circ\Phi)(\ox)\xi\big)$ and $\eta'\in N_{\ss{{\cal D}^*}}\big(\nabla(h\circ\Phi)(\ox)\xi\big)$ such that $\tet=\eta+\eta'$. Since ${\cal D}$ is a closed convex cone, we get $\eta'\in({\cal D}^{*})^*={\cal D}$ and hence $\eta'\in\ker\nabla(h\circ\Phi)(\ox)^*$. It follows from \eqref{io01} that $\eta\in DN_C(h(\oz),\bar\mu)(\nabla(h\circ\Phi)(\ox)\xi)$. Employing this together with \eqref{up10}, we arrive at the relationships
$$
\nabla_{x}\Psi^r(\ox,\omu)\xi+\nabla(h\circ\Phi)(\ox)^*\eta=0\;\mbox{ and }\;\eta\in DN_C\big(h(\oz),\bar\mu\big)(\nabla(h\circ\Phi)(\ox)\xi)\;\mbox{ with }\;\xi\ne 0,
$$
which contradict the noncriticality of $\omu$ and hence verifies \eqref{upper2}.

To finalized the proof, take the obtained constant $\ve'$ and the neighborhoods $V$ and $W$ from the Claim above and suppose without loss of generality that $\ve'<\al/2$ with $\al$ taken from \eqref{lm13}. Observe that there is a constant $\kappa\ge 0$ such that for any $(v,w)\in V\times W$ and any $(x_{vw},\mu_{vw})\in S^r(v,w)\cap\B_{\ve'}(\ox,\omu)$ we have the estimate
\begin{equation}\label{up2}
d\big(\mu_{vw};\Lm^r(\ox)\big)\le\kappa\big(\|x_{vw}-\ox\|+\|v\|+\|w\|\big).
\end{equation}
Indeed, \eqref{up2} can be justified by the same arguments as \eqref{lm14}. Combining \eqref{up2} and \eqref{upper2} gives us (\ref{rupper}) and thus verifies that the mapping $S^r$ from \eqref{rmapS} is semi-isolatedly calm at $\big((0,0),(\ox,\omu)\big)$. Invoking Lemma~\ref{axi3} tells that the semi-isolated calmness of the mapping $S^r$ yields the one for the mapping $S$ from \eqref{mapS}. This completes the proof of the theorem.
\end{proof}\vspace*{-0.07in}

Next we provide detailed discussions of our main result, Theorem~\ref{uplip}, and its proof.\vspace*{-0.05in}

\begin{Remark}[\bf discussing the obtained characterizations of noncriticality]\label{dis1} {\rm Our approach to characterize noncriticality of Lagrange multipliers for general variational systems \eqref{VS} developed above largely departs from those used in \cite[Theorem~1.43]{is14} and \cite[Theorem~4.1]{ms17} in polyhedral settings. Indeed, the proof of implication
(ii)$\Longrightarrow$(i) in Theorem~\ref{uplip} is significantly simplified due to the better translation of noncriticality via implication \eqref{cri2} that holds for any closed set $\Th$. The proof of (i)$\Longrightarrow$(ii) starts with a similar device as in the polyhedral case but departs from the latter in several steps. A new idea here is to deal with $\mu_k-\mu'_k$ instead of $\mu_k-\bar\mu$ to bypass the nonpolyhedrality of $\Th$. The term $\mu_k-\bar\mu$ works well in the proofs of \cite[Theorem~1.43]{is14} and \cite[Theorem~4.1]{ms17} due to  intrinsic properties of convex polyhedra, while using the same idea in nonpolyhedral cases of \cite[Theorem~3.3]{zz} and \cite[Proposition~4.2]{lp18} requires imposing strong assumptions, which may not hold even for the polyhedral settings of\cite{is14,ms17}. Our new proof of (i)$\Longrightarrow$(ii) resolves this issue by considering $\mu_k-\mu'_k$ and appealing to calculus of normal cones for convex cones under weak assumptions that holds in our setting due to the closedness assumption \eqref{lm15}. In this way a new term appears in our proof; namely,
\begin{equation}\label{lm17}
K_\Th(\oz,\olm)^*\cap\ker\nabla\Phi(\ox)^*,
\end{equation}
which is equivalent to $DN_\Th(\oz,\olm)(0)\cap\ker\nabla\Phi(\ox)^*$ due to the calculation of the graphical derivative of the normal cone mapping taken from Proposition~\ref{gdno}. As follows from  Theorem~\ref{unique}, this condition relates to uniqueness of the Lagrange multipliers. It appears naturally in our analysis and allows us to address generalized KKT systems with nonunique multipliers.

Observe further that the closedness assumption \eqref{lm15} is automatic if the set of Lagrange multipliers is a singleton and the mapping $M_{\ox}$ is calm at $((0,0),\olm)$. In this case we get from Theorem~\ref{unique} that the set in \eqref{lm17} is $\{0\}$, and thus \eqref{lm15} reduces to the closed set $K_\Th(\oz,\olm)^*$. Another important case where the assumed closedness holds is when $\Th$ is a convex polyhedron, which ensures the polyhedrality and hence closedness of $K_\Th(\oz,\olm)^*$. It is currently unclear whether the closedness of \eqref{lm15} is essential for the validity of (i)$\Longrightarrow$(ii) in Theorem~\ref{uplip}.

Note also that the calmness of the Lagrange multiplier mapping $M_{\ox}$ at $((0,0),\olm)$ assumed in Theorem~\ref{uplip}(b) always holds when
$\Th$ is a convex polyhedron. This condition is equivalent to the validity of \eqref{gf01} being a consequence of the Hoffman lemma; cf.\ Remark~\ref{rem1}. The following example shows that the calmness assumption on $M_{\ox}$ cannot be dropped in nonpolyhedral settings even in the case of unique Lagrange multipliers.}
\end{Remark}\vspace*{-0.13in}

\begin{Example}[\bf failure of noncriticality in the absence of calmness of Lagrange multipliers]\label{ex3}{\rm
Consider the semidefinite problem \eqref{sdp1} and recall from Example~\ref{ex1} that
$\Lambda_{\ss{c}}(\ox)=\{\olm\}$. It follows from Example~\ref{ex2} that the Lagrange multiplier mapping $M_{\ox}$ is not
calm at $((0,0),\olm)$. Further, we can conclude from \eqref{db01} that
\begin{eqnarray*}
&&K_{{\cal S}_+^2}(\oz,\olm)^*-K_{{\cal S}_+^2}(\oz,\olm)^*\cap\ker\nabla\Phi(\ox)^*\\
&=&\Big\{\begin{pmatrix}
b_{11}&b_{12}\\
b_{12}&b_{22}\\
\end{pmatrix}\in{\cal S}^2\;\Big|\;b_{22}\le 0\Big\}-\Big\{\begin{pmatrix}
0&a_{12}\\
a_{12}&0\\
\end{pmatrix}\in {\cal S}^2\;\Big|\;a_{12}\in\R\Big\}=\Big\{\begin{pmatrix}
b_{11}&b_{12}\\
b_{12}&b_{22}\\
\end{pmatrix}\in{\cal S}^2\;\Big|\;b_{22}\le 0\Big\},
\end{eqnarray*}
which ensures that the closedness assumption \eqref{lm15} of Theorem~\ref{uplip} is satisfied. Moreover, we know from Example~\ref{ex1} that
the unique Lagrange multiplier $\olm$ is noncritical. Our major goal is to show that the mapping $S$ from \eqref{mapS} for this problem is not semi-isolatedly calm at $\big((0,0),(\ox,\olm)\big)$, which demonstrates therefore that characterization (ii) of noncriticality of Lagrange multipliers in Theorem~\ref{uplip} fails without the calmness assumption on $M_{\ox}$. Observing that in the SDP framework \eqref{sdp1} the solution map $S$ reads as
\begin{eqnarray*}
S(v,w)=\big\{(x,\lm)\in\R^2\times{\cal S}^2\;\big|\;v=\nabla_x L(x,\lm),\;\lm\in N_{{\cal S}_+^2}(\Phi(x)+w)\big\}
\end{eqnarray*}
with $(v,w)\in\R^2\times{\cal S}^2$, we will actually get more: for any arbitrary small $t>0$ there are $(v_t,w_t)\in\B_t(0,0)\subset\R^2\times{\cal S}^2$ and $(x_t,\lm_t)\in S(v_t,w_t)\cap\B_t(\ox,\olm)$ such that {\em both} terms $\|\lm_t-\olm\|$ and $\|x_t-\ox\|$ are not of order $O(\|v_t\|+\|w_t\|)$; each of these properties yields the failure of the semi-isolated calmness of $S$ at $\big((0,0),(\ox,\olm)\big)$.

Considering first the $\lambda$-term, denote $v_t:=(-\frac{t^2}{2},-\frac{t^2}{2})$, $w_t:=\diag(0,0)$, $x_t:=\ox$, and
$\lm_t:=\begin{pmatrix}
-1-\frac{t^2}{2}& \frac{t}{2}\\
\frac{t}{2}&-\frac{t^2}{2}\\
\end{pmatrix}$
in the framework of Example~\ref{ex2}. As demonstrated therein, we have $\|\lm_t-\olm\|=O(t)$ while $O(\|v_t\|+\|w_t\|)=O(t^2)$. This verifies the claimed assertion on $\|\lm_t-\olm\|$ and confirms the failure of the semi-isolated calmness property for $S$ at $\big((0,0),(\ox,\olm)\big)$.

Next we show that the term $\|x_t-\ox\|$ also cannot be of order $O(\|v_t\|+\|w_t\|)$ in the absence of the calmness of the multiplier mapping $M_{\ox}$. This fact is instructive to understand the importance of the latter calmness property for superlinear convergence of primal iterations of SQP and related algorithms for solving nonpolyhedral conic programs. To proceed, denote $v_t:=(0,0)$ and
$w_t:=\begin{pmatrix}
0&t^2\\
t^2&0\\
\end{pmatrix}$
for which $O(\|v_t\| +\|w_t\|)=O(t^2)$ and then observe that $S$ can be considered as the KKT system for the parameterized semidefinite problem $P(t)$ given by
\begin{equation}\label{sdp2}
P(t):\quad\mbox{minimize}\;x_1+\frac{1}{2}x_1^2+\frac{1}{2}x_2^2\;\mbox{ subject to }\;\Phi(x)+w_t\in{\cal S}^2_+.
\end{equation}

It is proved in \cite[Example~4.5]{sh} (see also \cite[Example~4.54]{bs}) that the optimal solution mapping for \eqref{sdp2} is not outer Lipschitzian. Now we are going to verify the failure of the essentially more delicate semi-isolated calmness property of the solution map $S$ meaning that for the above pair $(v_t,w_t)$ there exists $(x_t,\lm_t)\in S(v_t,w_t)\cap\B_t(\ox,\olm)$ whenever $t>0$ is small enough. The latter task requires a significantly more involved analysis in comparison with \cite{sh}. We provide it below along with the verifying the aforementioned growth condition for $\|x_t-\ox\|$.

First observe that the parametric optimization problem \eqref{sdp2} is equivalent to
\begin{equation*}
\mbox{minimize}\;\ph_t(x_1,x_2):=\frac{1}{2}\Big((x_1+1)^2+x_2^2-1\Big)+\dd_{{\cal S}^2_+}\big(\Phi(x)+w_t\big)\;
\mbox{ subject to }\;x=(x_1,x_2)\in\R^2.
\end{equation*}
It is easy to see that the level sets of $\ph_t$ are uniformly bounded, which ensures the existence of minimizers for \eqref{sdp2} by the parametric version of the Weierstrass theorem; see, e.g., \cite[Theorem~1.17(a)]{rw}. Denote by $x_t=(x_{t1},x_{2t})$ such a minimizer for $P(t)$ and notice that the family $\{x_t\}$ as $t>0$ is uniformly bounded due to this property for the level sets of $\ph_t$.

Recall from Example~\ref{ex1} that $\ox$ is a unique minimizer for $P(0)$. Furthermore, it is clear from \eqref{lin1} that $(t^2,t^2)$ is a feasible solution to $P(t)$, and so
$$
0\le x_{t1}^2\le t^4+t^2\;\mbox{ and }\;0\le x_{t2}^2\le t^4+t^2,
$$
which yields $x_t\to \ox$ as $t\dn 0$. Note that the {\em Robinson constraint qualification} (RCQ)
$$
N_{{\cal S}^2_+}\big(\Phi(\ox)\big)\cap\ker\nabla\Phi(\ox)^*=\{0\}
$$
is satisfied for $P(0)$ and hence for $P(t)$ with small $t$ due to robustness of RCQ. This ensures that the set of Lagrange multipliers for $P(t)$ associated with $x_t$ is nonempty and uniformly bounded if $t$ is sufficiently small. Thus there is $\ve>0$ and $l\ge 0$ with
\begin{equation}\label{lin0}
\|\lm_t\|\le l\;\mbox{ whenever }\;|t|\le\ve
\end{equation}
for such Lagrange multipliers. It follows from $\Lambda_{\ss{c}}(\ox)=\{\olm\}$ that $\lm_t\to\olm$ as $t\dn 0$ and
$(x_t,\lm_t)\in S(v_t,w_t)$. Letting $\lm_t:=\begin{pmatrix}
\lm^t_{11}&\lm^t_{12}\\
\lm^t_{12}&\lm^t_{22}
\end{pmatrix}$, obtain from the first-order optimality conditions that
\begin{equation}\label{lin00}
v_t=\nabla_x L(x_t,\lm_t)\iff\lm^t_{11}=-x_{t1}-1,\;\lm^t_{22}=-x_{t2}\;\mbox{ and}
\end{equation}
\begin{equation*}
\lm_t\in N_{{\cal S}_+^2}(\Phi(x_t)+w_t)\iff\Phi(x_t)+w_t\in{\cal S}_+^2,\;\lm_t\in{\cal S}_-^2,\;\lm_t\big(\Phi(x_t)+w_t\big)=\diag(0,0).
\end{equation*}
The latter tells us by elementary linear algebra that
\begin{equation}\label{lin1}
\Phi(x_t)+w_t=\begin{pmatrix}
x_{t1}&t^2\\
t^2&x_{t2}\\
\end{pmatrix}\in{\cal S}_+^2\iff x_{t1}\ge 0,\;x_{t2}\ge 0,\;x_{t1}x_{t2}\ge t^4\;\mbox{ and}
\end{equation}
\begin{equation}\label{lin2}
\lm_t=\begin{pmatrix}
\lm^t_{11}&\lm^t_{12}\\
\lm^t_{12}&\lm^t_{22}\\
\end{pmatrix}\in{\cal S}_-^2\iff\lm^t_{11}\le 0,\;\lm^t_{22}\le 0,\;\lm^t_{11}\lm^t_{22}\ge 3(\lm^t_{12})^2.
\end{equation}
Moreover, it follows from $\lm_t(\Phi(x_t)+w_t)=\diag(0,0)$ that
\begin{equation}\label{lin3}
\lm^t_{11}x_{t1}+t^2\lm^t_{12}=0,\;\lm^t_{22}x_{t2}+t^2\lm^t_{12}=0,\;t^2\lm^t_{11}+x_{t2}\lm^t_{12}=0,\;\mbox{ and }\;t^2\lm^t_{22} +x_{t1}\lm^t_{12}=0.
\end{equation}
Using the first two equations in \eqref{lin3} together with \eqref{lin1} and \eqref{lin2} implies that
\begin{equation}\label{lin4}
x_{t1}x_{t2}=t^4\;\mbox{ and }\;\lm^t_{11}\lm^t_{22}=(\lm^t_{12})^2.
\end{equation}
The latter tells us, being combined with the last two equations in \eqref{lin3}, that $\lm^t_{22}x_{t2}=\lm^t_{11}x_{t1}$, which yields in turn the relationship
\begin{equation}\label{lin5}
x_{t2}^3=-\lm^t_{11}t^4.
\end{equation}
This along with \eqref{lin0} verifies that $|x_{t2}|=O(t^{\frac{4}{3}})$ and hence allows us to deduce from $\lm_t\to\olm$ as $t\dn 0$ that $|\lm^t_{11}|\ge\frac{1}{2}$ for all $t$ sufficiently small. Using it and the first equation in \eqref{lin4} together with \eqref{lin5}, we get $x_{t1}\sqrt[3]{-\lm^t_{11}}=t^{\frac{8}{3}}$ and so arrive at $|x_{t1}|=O(t^{\frac{8}{3}})$. Employing the latter condition together with \eqref{lin5} again brings us to
$$
\|x_t-\ox\|=\|x_t\|=O(t^{\frac{4}{3}})\;\mbox{ and }\;x_t\in\B_{t/2}(\ox)\;\mbox{ for all small }\;t>0.
$$
Combining it with \eqref{lin00} and the second equation in \eqref{lin4} shows that
$$
\|\lm_t-\olm\|=O(t^{\frac{4}{3}})\;\mbox{ and }\;\lm_t\in\B_{\frac{t}{2}}(\olm)\;\mbox{ whenever }\;t\;\mbox{ is sufficiently small}.
$$
This tells us that $(x_t,\lm_t)\in S(v_t,w_t)\cap\B_t(\ox,\olm)$, that both terms $\|x_t-\ox\|$ and $\|\lm_t-\olm\|$ are of order $O(t^{\frac{4}{3}})$, and therefore
$$
\lim_{t\dn 0}\frac{\|x_t-\ox\|}{\|v_t\|+\|w_t\|}=\lim_{t\dn 0}\frac{\|\lm_t-\olm\|}{\|v_t\|+\|w_t\|}=\lim_{t\dn 0}\frac{O(t^{\frac{4}{3}})}{\sqrt{2}\,t^2}=\infty.
$$
It verifies all the claims made above and thus confirms that the calmness of the Lagrange multiplier mapping is essential for the obtained characterizations of noncritical multipliers in nonpolyhedral variational systems.}
\end{Example}\vspace*{-0.05in}

The next result strongly relates to Theorem~\ref{uplip} while giving us a significant additional information. It shows that a new second-order condition, which strengthens noncriticality, yields the semi-isolated calmness property of the solution map \eqref{mapS} at $((0,0),(\ox,\olm))$ without imposing the closedness assumption while providing that the multiplier mappings $M_{\ox}$ is calm at $((0,0),\olm)$. The new {\em second-order condition} for \eqref{VS} reads as follows:
\begin{equation}\label{nsoc}
\begin{cases}
\big\la\nabla_{x}\Psi(\ox,\olm)\xi,\xi\big\ra+\big\la\nabla^2\la\bar\mu,h\ra(\oz\nabla\Phi(\ox)\xi,\nabla\Phi(\ox)\xi\big\ra>0 \\
\mbox{for all }\;0\ne\xi\in\X\;\mbox{ with }\;\nabla\Phi(\ox)\xi\in K_\Th(\oz,\olm),
\end{cases}
\end{equation}
where $h$ and $\omu$ are taken from \eqref{red} and \eqref{redl}, respectively. When $\Phi=\nabla_x L$ with $L$ standing for the standard Lagrangian in constrained optimization \eqref{coop}, condition \eqref{nsoc} reduced to the second-order sufficient condition \eqref{ssoc}.\vspace*{-0.05in}

\begin{Theorem}[\bf semi-isolated calmness from second-order condition]\label{socic} Let $(\ox,\olm)$ be a solution to \eqref{VS}, let
$\Th$ be ${\cal C}^2$-cone reducible at $\oz=\Phi(\ox)$ to a closed convex cone $C$, and let the multiplier mapping $M_{\ox}$ from \eqref{lagmap} be calm at $((0,0),\olm)$. If the second-order condition \eqref{nsoc} holds, then the solution map $S$ from \eqref{mapS} is semi-isolatedly calm at $((0,0),(\ox,\olm))$.
\end{Theorem}\vspace*{-0.18in}
\begin{proof} We utilize a reduction procedure similar to the device of Theorem~\ref{uplip} and thus present just a sketch of the proof. Considering the reduced system \eqref{rVS}, observe that \eqref{nsoc} corresponds to the reduced second-order condition
\begin{equation}\label{rnsoc}
\big\la\nabla_{x}\Psi^r(\ox,\omu)\xi,\xi\big\ra>0\;\mbox{ for all }\;0\ne\xi\in\X\;\mbox{ with }\;\nabla(h\circ \Phi)(\ox)\xi\in K_C\big(h(\oz),\omu\big)
\end{equation}
for \eqref{rVS}; see \cite[equation~(3.272)]{bs} for more detail. By Lemma~\ref{axi4} it suffices to show that the solution map $S^r$ from \eqref{rmapS} is semi-isolated calm at $((0,0),(\ox,\omu))$. To this end, we proceed as the proof of Theorem~\ref{uplip} and show first that \eqref{upper2} fulfills. Arguing by contradiction and proceeding as in the proof of Theorem~\ref{uplip} give us
\eqref{up10}, \eqref{tr04}, and \eqref{mn02} without using the closedness condition \eqref{lm15}. It implies in turn that
$$
0=\la 0,\xi\ra=\la\nabla_{x}\Psi^r(\ox,\omu)\xi,\xi\ra+\la\tet,\nabla(h\circ\Phi)(\ox)\xi\ra\ge\la\nabla_{x}\Psi^r(\ox,\omu)\xi,\xi\ra
$$
with $\xi\ne 0$ due to \eqref{xi} and $\tet$ taken from \eqref{seq1}. Employing \eqref{tr04} along with \eqref{rnsoc} yields $\xi=0$, a contradiction, which verifies \eqref{upper2}. Finally, we can justify \eqref{up2} as in the proof of Theorem~\ref{uplip} using the calmness of the multiplier mapping $M_{\ox}$ at $((0,0),\olm)$.
\end{proof}\vspace*{-0.1in}
 
In the constrained optimization framework \eqref{coop}, the obtained result provides an important extension of the fact well recognized for NLPs. Indeed, it can be distilled from \cite[Lemma~2]{hg} that the second-order sufficient condition \eqref{ssoc} yields the semi-isolated
calmness of $S$. Theorem~\ref{socic} reveals that such a result can be guaranteed in the general framework of \eqref{VS} if in addition to the second-order condition \eqref{rnsoc} the Lagrange multiplier mapping $M_{\ox}$ is calm. Remember that the latter property is automatic for NLPs.
Moreover, combining Examples~\ref{ex1} and \ref{ex3} tells us that the calmness of $M_{\ox}$ is essential in Theorem~\ref{socic}.\vspace*{0.02in}

The final result of this section provides an efficient condition ensuring the validity of both assumptions on closedness \eqref{lm15} and calmness of Lagrange multipliers imposed in Theorem~\ref{uplip}(b). In this way we get complete characterizations of noncriticality of Lagrange multipliers via the error bound and semi-isolated calmness of solution maps to nonpolyhedral systems as in the case of polyhedrality. The condition we are going to use is known as {\em strict complementarity} \cite[Definition~4.74]{bs} for \eqref{VS} at $\ox$ meaning that there is $\lm\in\Lambda(\ox)$ such that $\lm\in\ri N_\Th(\Phi(\ox))$.\vspace*{-0.1in}

\begin{Theorem}[\bf characterizations of noncriticality of multipliers under strict complementarity]\label{suffc} Let $\ox$ be a stationary point from \eqref{stat}, let $\Th$ be ${\cal C}^2$-cone reducible at $\oz=\Phi(\ox)$ to a closed convex cone $C$, and let the strict complementarity condition be satisfied at $\ox$ for \eqref{VS} at $\ox$. Then a Lagrange multiplier $\olm\in\Lambda(\ox)$ is noncritical if and only if either one of the conditions {\rm(ii)} and {\rm(iii)} of Theorem~{\rm\ref{uplip}} is satisfied.
\end{Theorem}\vspace*{-0.18in}
\begin{proof} This theorem follows from Theorem~\ref{uplip} provided that the imposed strict complementarity implies both the closedness condition \eqref{lm15} and the calmness of the multiplier mapping $M_{\ox}$ assumed in Theorem~\ref{uplip}(b). We split the proof into the following three steps.\\[1ex]
{\bf Step~1:} {\em The strict  complementarity condition holds for \eqref{VS} if and only if it holds for the reduced KKT system \eqref{rVS}.} To verify this claim, suppose that the strict complementarity condition holds at $\ox$ for \eqref{VS} and then find a multiplier $\lm\in\Lambda(\ox)$ such that $\lm\in\ri N_\Th(\Phi(\ox))$. It follows from the normal cone calculus \eqref{rch} and from \cite[Proposition~2.44]{rw} that
\begin{equation*}
\lm\in\ri N_\Th(\oz)=\ri\big(\nabla h(\oz)^*N_C(h(\oz))\big)=\nabla h(\oz)^*\big(\ri N_C(h(\oz))\big).
\end{equation*}
This ensures the existence of a vector $\mu\in\ri N_C(h(\oz))$ such that $\lm=\nabla h(\oz)^*\mu$. Unifying this with $\lm\in\Lambda(\ox)$ gives us $\mu\in\Lambda^r(\ox)$ and shows therefore that the strict complementarity condition holds for \eqref{rVS}. The opposite implication is proved similarly.\\[1ex]
{\bf Step~2:} {\em The strict complementarity condition for \eqref{VS} at $\ox$ yields the closedness condition in Theorem~{\rm\ref{uplip}(b)}}. It follows from Step~1 that we need to verify the closedness of the set
\begin{equation}\label{hg01}
K_C\big(h(\oz),\bar\mu\big)^*-\big[K_C\big(h(\oz),\bar\mu\big)^*\cap\ker\nabla(h\circ\Phi)(\ox)^*\big]
\end{equation}
from Lemma~\ref{axi1}(ii) under the validity of the strict complementarity condition for the reduced system \eqref{rVS}. To furnish this, recall that $h(\oz)=0$, and hence $K_C(h(\oz),\bar\mu)^*=\cl(C^*+\R\omu)$. Since $\omu\in C^*$ and $\span C^*$ is closed, we have $\cl (C^*+\R\omu)\subset\span C^*$. This leads us to
\begin{eqnarray*}
K_C\big(h(\oz),\bar\mu\big)^*-\big[K_C\big(h(\oz),\bar\mu\big)^*\cap\ker\nabla(h\circ\Phi)(\ox)^*\big]&\subset&K_C\big(h(\oz),\bar\mu\big)^*-
K_C\big(h(\oz),\bar\mu\big)^*\\
&=&\cl(C^*+\R\omu)-\cl(C^*+\R\omu)\\
&\subset&\span C^*-\span C^*=\span C^*.
\end{eqnarray*}
On the other hand, it follows from the strict complementarity condition for \eqref{rVS} that there is a vector $\mu\in\ri N_C(h(\oz))=\ri C^*$ such that $\mu\in\Lambda^r(\ox)$. Pick $w\in\span C^*$ and observe that $\aff C^*=\span C^*$. By $\mu\in\ri C^*$ we find a small number
$t>0$ for which $\mu+tw\in C^*$. Combining the above facts brings us to the relationships
\begin{eqnarray*}
tw=(\mu+t w-\omu)-(\mu-\omu)&\subset&\cl(C^*+\R\omu)-\big[\cl(C^*+\R\omu)\cap\ker\nabla(h\circ\Phi)(\ox)^*\big]\\
&\subset&K_C\big(h(\oz),\bar\mu\big))^*-\big[K_C\big(h(\oz),\bar\mu\big)^*\cap\ker\nabla(h\circ\Phi)(\ox)^*\big],
\end{eqnarray*}
which readily imply the inclusion
\begin{equation*}
\span C^*\subset K_C\big(h(\oz),\bar\mu\big)^*-\big[K_C\big(h(\oz),\bar\mu\big)^*\cap\ker\nabla(h\circ\Phi)(\ox)^*\big].
\end{equation*}
Since the opposite inclusion also holds by the above discussion, we come up with the equality
\begin{equation*}
K_C\big(h(\oz),\bar\mu\big)^*-\big[K_C\big(h(\oz),\bar\mu\big)^*\cap\ker\nabla(h\circ\Phi)(\ox)^*\big]=\span C^*,
\end{equation*}
which verifies the closedness of the set in \eqref{hg01}. Appealing now to Lemma~\ref{axi1} tells us that the set in \eqref{lm15} is closed as well.\\[1ex]
{\bf Step~3:} {\em The strict complementarity condition for \eqref{VS} at $\ox$ implies that the multiplier set $M_{\ox}$ is calm at $((0,0),\olm)$}. By Step~1 it suffices to prove that estimate \eqref{lm13} holds under the strict complementarity condition for \eqref{rVS}. Remembering that $h(\oz)=0$ gives us $\Lambda^r(\ox)=\{\mu\in\E\;|\;\Psi^r(\ox,\mu)=0,\;\mu\in C^*\}$. This together with \cite[Corollary~3]{bbw} and the Hoffman lemma ensures the existence of numbers $\ve>0$ and $\ell\ge 0$ for which
\begin{equation}\label{fm00}
d\big(\mu;\Lambda^r(\ox)\big)\le\ell\big(\|\Psi^r(\ox,\mu)\|+d(\mu;C^*)\big)\;\mbox{ whenever }\;\mu\in\B_{\ve}(\omu).
\end{equation}
Pick $\mu\in \E$ and let $y:=P_C(\mu)$, where $P_C(\mu)$ stands for the projection of $\mu$ onto the convex cone $C$. It implies that $\mu-y\in N_C(y)$ and so $\mu-y\in C^*$, which brings us to
\begin{equation}\label{fm01}
d(\mu;C^*)\le\|\mu-(\mu-y)\|=\|y\|=\|P_C(\mu)\|\;\mbox{ for all }\;\mu\in\E.
\end{equation}
On the other hand, we get that $P_C(\mu)=0$ if and only if $\mu\in C^*$. This allows us to deduce from $\mu\in C^*$ the equalities
\begin{equation*}
\|P_C(\mu)\|=0=d\big(0;N_{C^*}(\mu)\big)=d\big(0;N^{-1}_{C}(\mu)\big)=d\big(h(\oz);N^{-1}_{C}(\mu)\big).
\end{equation*}
If $\mu\not\in C^*$, then $\|P_C(\mu)\|<d(h(\oz);N_{C^*}(\mu))=d(h(\oz);N^{-1}_{C}(\mu))=\infty$, and so
\begin{equation}\label{fm02}
\|P_C(\mu)\|\le d\big(h(\oz);N^{-1}_{C}(\mu)\big)\;\mbox{ for all }\;\mu\in\E.
\end{equation}
Combining \eqref{fm00}--\eqref{fm02} verifies estimate \eqref{lm13}, which yields by Lemma~\ref{axi3} the calmness of the Lagrange multiplier mapping $M_{\ox}$ at $((0,0),\olm)$ and thus completes the proof.
\end{proof}\vspace*{-0.1in}

It follows from \cite{sh97} that the strict complementarity condition ensures the equivalence between the uniqueness of Lagrange multipliers and the strong Robinson constraint qualification \eqref{sqc} for problems of semidefinite programming. Theorem~\ref{suffc} allows us to extend Shapiro's result to the general ${\cal C}^2$-cone reducible setting of \eqref{VS}.\vspace*{-0.1in}

\begin{Corollary}[\bf uniqueness of Lagrange multipliers under the strict complementarity condition]\label{shap} Let $(\ox,\olm)$ be a solution to the variational system \eqref{VS}, where $\Th$ is ${\cal C}^2$-cone reducible at $\oz=\Phi(\ox)$ to a closed convex cone $C$. Assume that
the strict complementarity condition holds at $\ox$ for \eqref{VS}. Then the Lagrange multiplier set $\Lambda(\ox)$ is a singleton if and only if the equivalent qualification conditions \eqref {gf02} and \eqref{sqc} are satisfied.
\end{Corollary}\vspace*{-0.18in}
\begin{proof} This follows from the combination of Theorems~\ref{unique}, \ref{suffc} and Proposition~\ref{pdcq}.
\end{proof}
\end{document}